%% file: article.tex
\documentclass{gtart}

\newcommand{\pathtotrunk}{./}
\input{text/article_preamble.tex}
\input{text/top_matter.tex}

\begin{document}

\begin{abstract}
We give a combinatorial description of the ``$D_{2n}$ planar algebra,''
by generators and relations. We explain how the generator
interacts with the Temperley-Lieb braiding.
This shows the previously known braiding on the even part
extends to a `braiding up to sign' on the entire planar algebra.

We give a direct proof that our relations are consistent (using this
`braiding up to sign'), give a complete description of the associated
tensor category and principal graph, and show that the planar algebra
is positive definite. These facts allow us to identify our combinatorial
construction with the standard invariant of the subfactor $D_{2n}$.
\end{abstract}

\maketitle

\setcounter{tocdepth}{2}%
\tableofcontents

\section{Introduction}\label{sec:intro}
\input{text/intro.tex}

\section{Background}
\input{text/background}

\section{Skein theory}
\input{text/presentation.tex}
\input{text/braiding.tex}
\input{text/consistency.tex}

\section{The planar algebra $\pa$ is $D_{2n}$}
This planar algebra is called $D_{2n}$ because it is the unique subfactor
planar algebra with principal graph $D_{2n}$.  To prove this we will need
two key facts: first that its principal graph is $D_{2n}$, and second that
it is a subfactor planar algebra.

In \S \ref{sec:category}, we describe a tensor category associated to any
planar algebra, and using that define the principal graph. We then check
that the principal graph for $\pa$ is indeed the Dynkin
diagram $D_{2n}$.

In \S \ref{sec:basis}, we exhibit an explicit basis for the planar
algebra.  This makes checking positivity straightforward.

\input{text/category.tex}
\input{text/basis.tex}

\input{text/appendix.tex}

\newcommand{\urlprefix}{}
\bibliographystyle{unsrt}
\bibliography{bibliography/bibliography}

This paper is available online at \arxiv{0808.0764}, and at
\url{http://tqft.net/d2n}.

\end{document}

%% file: text/article_preamble.tex
\input{\pathtotrunk preamble.tex}

\ifpdf
\usepackage[pdftex]{graphicx}
\else
\usepackage[dvips]{graphicx}
\fi

\usepackage{tikz}
\usetikzlibrary{shapes}

\usepackage{color}

%% file: preamble.tex

\usepackage{amsmath,amssymb,amsfonts}
\usepackage{ifpdf}
\usepackage[usenames]{color}
\usepackage{enumerate}

\usepackage{leftidx}

\ifpdf
\usepackage[pdftex,all,color]{xy}
\else
\usepackage[all,color]{xy}
\fi

\SelectTips{cm}{}

\vfuzz2pt 
\hfuzz2pt 

\def\RCS$#1: #2 ${\expandafter\def\csname RCS#1\endcsname{#2}}
\RCS$Revision: 543 $ \RCS$Date: 2009-02-07 08:31:38 +1100 (Sat, 07 Feb 2009) $

\ifpdf
\newcommand{\pathtodiagrams}{\pathtotrunk diagrams/pdf/}
\else
\newcommand{\pathtodiagrams}{\pathtotrunk diagrams/eps/}
\fi

\newcommand{\mathfig}[2]{{\hspace{-3pt}\begin{array}{c}%
  \raisebox{-2.5pt}{\includegraphics[width=#1\textwidth]{\pathtodiagrams #2}}%
\end{array}\hspace{-3pt}}}

\newcommand{\rotatemathfig}[3]{{\hspace{-3pt}\begin{array}{c}%
  \raisebox{-2.5pt}{\rotatebox{#2}{\includegraphics[height=#1\textwidth]{\pathtodiagrams #3}}}%
\end{array}\hspace{-3pt}}}

\ifpdf
\usepackage[pdftex,plainpages=false,hypertexnames=false,pdfpagelabels]{hyperref}
\else
\usepackage[dvips,plainpages=false,hypertexnames=false]{hyperref}
\fi
\newcommand{\arxiv}[1]{\href{http://arxiv.org/abs/#1}{\tt arXiv:\nolinkurl{#1}}}
\newcommand{\doi}[1]{\href{http://dx.doi.org/#1}{{\tt DOI:#1}}}

\newcommand{\googlebooks}[1]{(preview at \href{http://books.google.com/books?id=#1}{google books})}
\newcommand{\numdam}[1]{}

\theoremstyle{plain}
\newtheorem*{fact}{Fact}
\newtheorem{prop}{Proposition}[section]

\newtheorem{thm}[prop]{Theorem}
\newtheorem{lem}[prop]{Lemma}
\newtheorem{cor}[prop]{Corollary}
\newtheorem*{cor*}{Corollary}

\newtheorem{defn}[prop]{Definition}         
\newtheorem*{defn*}{Definition}             

\newtheorem*{kuperberg}{The Kuperberg Program}
\newtheorem*{mainthm}{The Main Theorem}

\newenvironment{rem}{\noindent\textsl{Remark.}}{}  
\numberwithin{equation}{section}

\newcommand{\noqed}{\renewcommand{\qedsymbol}{}}

\newcounter{comment}
\newcommand{\noop}[1]{}

\def\clap#1{\hbox to 0pt{\hss#1\hss}}


\newcommand{\Complex}{\mathbb C}

\newcommand{\id}{\boldsymbol{1}}
\renewcommand{\imath}{\mathfrak{i}}
\renewcommand{\jmath}{\mathfrak{j}}

\newcommand{\iso}{\cong}

\newcommand{\abs}[1]{\left|#1\right|}

\newcommand{\card}[1]{\sharp{#1}}

\newcommand{\eset}{\emptyset}

\newcommand{\cC}{\mathcal{C}}

\newcommand{\qiq}[2]{[#1]_{#2}}
\newcommand{\qi}[1]{\qiq{#1}{q}}

\newcommand{\directSum}{\oplus}
\newcommand{\DirectSum}{\bigoplus}
\newcommand{\tensor}{\otimes}
\newcommand{\Tensor}{\bigotimes}

\newcommand{\Mat}[1]{\operatorname{\mathbf{Mat}}\left(#1\right)}

\newcommand{\Hom}[3]{\operatorname{Hom_{#1}}\left(#2,#3\right)}

\newcommand{\pa}{\mathcal{PA}(S)}
\newcommand{\TL}{\mathcal{TL}}
\newcommand{\JW}[1]{f^{(#1)}}
\newcommand{\tr}[1]{\text{tr}(#1)}


%% file: text/top_matter.tex
\title{Skein theory for the $D_{2n}$ planar algebras}

\author{Scott~Morrison}
\address{
   Microsoft Station Q, University of California, Santa Barbara 93106-6105
}%
\email{scott@tqft.net}

\author{Emily~Peters\thanks{Supported in part by NSF Grant DMS0401734}
}
\address{
   Department of Mathematics,
   University of California, Berkeley, 94720
}%
\email{eep@math.berkeley.edu}

\author{Noah~Snyder\thanks{Supported in part by RTG grant DMS-0354321}
}
\address{
   Department of Mathematics,
   University of California, Berkeley, 94720
}%
\email{nsnyder@math.berkeley.edu}
\urladdr{\url{http://tqft.net/} \url{http://math.berkeley.edu/~eep} \url{http://math.berkeley.edu/~nsnyder}}

\date{
  First edition: August 6 2008
}

\primaryclass{46L37} \secondaryclass{57M15 46M99} \keywords{
  Planar Algebras, Subfactors, Skein Theory
}

%% file: text/intro.tex
Start with a category with tensor products and a good theory of duals (technically a spherical tensor category \cite{MR1686423}, or slightly more generally a spherical $2$-category\footnote{Recall that a monoidal category is just a $2$-category with one object. Subfactor planar algebras are $2$-categories with two objects but still have a good spherical theory of duals for morphisms. In general one could consider any $2$-category with a good theory of duals. Watch out that our terminology differs from \cite{MR1686421}. There the phrase `spherical $2$-category' refers to a \emph{monoidal} $2$-category which we think would better be called a spherical monoidal $2$-category or a spherical $3$-category.}), such as the category of representations of a quantum group, or the
category of bimodules coming from a subfactor.  Fix your favorite object in this tensor category. Then the $\operatorname{Hom}$-spaces between arbitrary tensor products of the chosen object and its dual fit together into a structure called a planar algebra (a notion due to Jones \cite{math.QA/9909027}) or the roughly equivalent structure called a spider (a notion due to Kuperberg \cite{MR1403861}). 
 Encountering such an object should tempt you to participate in:

\begin{kuperberg}
Give a presentation by generators and relations for every interesting planar algebra.
Generally it's easy to guess some generators, and not too hard to
determine that certain relations hold. You should then aim to prove that
the combinatorial planar algebra given by these generators and relations
agrees with your original planar algebra.
Ideally, you also understand other properties of the original category (for example positivity, being spherical, or
being braided)
in terms of the presentation.
\end{kuperberg}

The difficulty with this approach is often in proving
combinatorially that your relations are self-consistent,
without appealing to the original planar algebra.
Going further, you could try to find
 explicit `diagrammatic' bases for all the original $\operatorname{Hom}$
spaces, as well as the combinatorial details of $6-j$ symbols or
`recombination' rules.

This program has been fulfilled completely for the $A_n$ subfactors (equivalently, for the representation theory of $U_q(\mathfrak{sl}_2)$ at a root of unity),
for all the subfactors coming from Hopf algebras \cite{MR2079886, MR2146224, MR1971553, MR2346882},
and for the representation categories of the rank $2$ quantum groups
\cite{MR1403861, MR2360947}. Some progress has been made on the representation categories of $U_q(\mathfrak{sl}_n)$ for $n \geq 4$ \cite{MR1659228, dongseok-thesis, scott-thesis}.
Other examples of planar algebras which have been described or constructed by generators and relations include  the BMW and Hecke algebras \cite{math.QA/9909027, MR1733737}, the Haagerup subfactor \cite{0902.1294}, and the Bisch-Haagerup subfactors \cite{MR1386923, 0807.4134}.

In this paper we apply the Kuperberg program to the subfactor planar
algebras corresponding to $D_{2n}$. The $D_{2n}$ subfactors are one of
the two infinite families (the other being $A_n$) of subfactors of index
less than $4$.  Also with index less than $4$ there are two sporadic
examples, the $E_6$ and $E_8$ subfactors.  See \cite{MR999799, MR996454,
MR1145672, MR1936496} for the story of this classification.

The reader familiar with quantum groups should be warned that although $D_{2n}$ is related to the Dynkin diagram $D_{2n}$, it is not in any way related to the quantum group $U_q(\mathfrak{so}_{4n})$.  To get from $U_q(\mathfrak{so}_{4n})$ to the $D_{2n}$ diagram you look at its roots. To get from the $D_{2n}$ subfactor to the $D_{2n}$ diagram you look at its fusion graph. The fusion graph of a quantum group is closely related to its fundamental alcove, not to its roots.  Nonetheless the $D_{2n}$ subfactor is related to quantum groups!  First, It is a quantum subgroup of $U_q(\mathfrak{sl}_2)$ in the sense of \cite{MR1936496}.  To make matters even more confusing, the $D_{2n}$ subfactor is related via level-rank duality to the quantum group $U_q(\mathfrak{so}_{2n-2})$; see \cite{d2n-links} for details.

The $D_{2n}$ subfactors were first constructed in
\cite{MR1308617}, using an automorphism of the subfactor $A_{4n-3}$.
(This `orbifold method' was studied further in
\cite{MR1379298,MR1309549}.) Since then, several papers have offered
alternative constructions; via planar algebras, in \cite{MR1929335}, and as
a module category over an algebra object in $A_{4n-3}$, in
\cite{MR1936496}. In this paper we'll show an explicit description of the
associated $D_{2n}$ planar algebra, and via the results of
\cite{math.QA/9909027,MR1334479} or of \cite{0807.4146} this gives an
indirect construction of the subfactor itself. 

Our goal in this paper is to understand as much as possible about the the $D_{2n}$ planar algebra on the level of planar algebras -- that is, without appealing to subfactors, or any structure beyond the combinatorics of diagrams.  
 We also hope
that our treatment of the planar algebra for $D_{2n}$ by generators and
relations nicely illustrates the goals of the Kuperberg program,
although more complicated examples will require different methods.

Our main object of study is a planar algebra $\pa$ defined by generators and relations.

\begin{defn}\label{def:pa}
Fix $q=\exp(\frac{\pi i}{4n-2})$.
Let $\pa$ be the planar algebra generated by a single ``box" $S$ with $4n-4$ strands, modulo the following relations.
\begin{enumerate}

\item\label{delta} A closed circle is equal to $[2]_q = (q+q^{-1}) = 2 \cos(\frac{\pi}{4n-2})$ times the empty diagram.

\item\label{rotateS} Rotation relation:
\input{diagrams/tikz/rotateS.tex}

\item\label{capS} Capping relation:
\input{diagrams/tikz/capS.tex}

\item\label{twoS} Two $S$ relation:
\input{diagrams/tikz/twoS.tex}

\end{enumerate}
\end{defn}

This paper uses direct calculations on diagrams  to establish the following theorem:

\begin{mainthm}
\label{thm:main}%
$\pa$ is the $D_{2n}$ subfactor planar algebra; that is,
\begin{enumerate}
\item the space of closed diagrams is $1$-dimensional,
\item $\pa$ is spherical,
\item the principal graph of $\pa$ is the Dynkin diagram $D_{2n}$, and
\item $\pa$ is unitary, that is, it has a star structure for which $S^*=S$, and the associated inner product is positive definite.
\end{enumerate}
\end{mainthm}

Many of the terms appearing in this statement will be given definitions later, although a reader already acquainted with the theory of subfactors should not
find anything unfamiliar.\footnote{Although perhaps they should watch out --- we'll define the principal graph of a planar algebra by a slightly
different route than usual, failing to mention either the basic construction \cite{MR1473221}, or bimodules over a von Neumann algebra \cite{MR1424954}!}

In this paper our approach is to start with the generators and relations for $\pa$ and to prove the Main Theorem from scratch.  The first part of the Main Theorem in fact comes in two subparts; first
that the relations given in
Definition \ref{def:pa} are consistent (that is, $\pa_0 \neq 0$),
 and second 
 that every closed diagram can be evaluated as a multiple of the empty
diagram using the relations.
 These
statements appear as Corollary
\ref{cor:evaluation} and as Theorem \ref{thm:consistency}. Corollary \ref{cor:spherical} proves that $\pa$ is
spherical.

Our main tool in showing all of this is a `braiding up to sign' on the
entire planar algebra $D_{2n}$; the details are in Theorem
\ref{thm:passacrossS}. It is well-known that the even part of $D_{2n}$ is
braided (for example \cite{MR1936496}), but we extend that braiding to
the whole planar algebra with the caveat that if you pull $S$ \emph{over} a
strand it becomes $-S$. In a second paper \cite{d2n-links}, we will give
results about the knot and link invariants which can be constructed using
this planar algebra. From these, we can derive a number of new identities
between classical knot polynomials.

In Section \ref{sec:category}, we will describe the structure of the tensor category of projections, essentially rephrasing the concepts of fusion algebras in planar algebra language.  Some easy diagrammatic calculations then establish the third part of the main theorem.   Section \ref{sec:basis} exhibits an orthogonal basis for the planar algebra, and the final part of the main theorem becomes an easy consequence.
{Finally, Appendix \ref{appendix} describes a family of related planar algebras, and sketches the corresponding results.}

In addition to our direct approach, one could also prove the main theorem in the following indirect way.  First take one of the known constructions of the subfactor $D_{2n}$.  By \cite{math.QA/9909027} the standard invariant of $D_{2n}$ gives a planar algebra. Using the techniques in \cite{MR1929335} and \cite{quadratic}, find the generator and some of the relations for this planar algebra.  At this point you'll have reconstructed our list of generators and relations for $\pa$.  However, even at this point you will only know that the $D_{2n}$ planar algebra is a quotient of $\pa$.  To prove that $D_{2n} = \pa$ you would still need many of the techniques from this paper.  In particular, using all the above results only allows you to skip Section 3.3 and parts of Section 4.2 (since positive definiteness would follow from non-degeneracy of the inner product and positivity for $D_{2n}$).

We'd like to thank Stephen Bigelow, Vaughan Jones, and Kevin Walker for interesting conversations.
During our work on this paper, Scott Morrison was at Microsoft Station Q, Emily Peters was supported in part by NSF Grant DMS0401734 and Noah Snyder was supported in part by RTG grant DMS-0354321.

We've fixed some errors pointed out by careful readers in earlier versions of this paper. Kevin Walker pointed out an error in the statement (but, happily, not the proof!) of Theorem \ref{easyconsequences}. Shohei Matsunaga pointed out an error in Equation \eqref{eq:wjwk-large}.
Kazuyo Sakamaki pointed out an error in the coefficients appearing in
Definition \ref{defn:algorithm} and Theorem
\ref{thm:algorithm-well-defined}.

%% file: diagrams/tikz/rotateS.tex
$
\begin{tikzpicture}[baseline]
    \clip (-.6,-1) rectangle (.6,1);

    \draw (-.32,0) -- ++(90:3mm) arc (0:180:1mm) -- ++(-90:3.5mm) arc (-180:0:5.3mm) -- ++(90:15mm);
    \draw (.32,0) -- ++(90:15mm);
    \draw (-.15,0) -- ++(90:15mm);
    \draw (.07,.7) node {...};
    
    \node (S) at (0,0) [circle, draw, fill=white] {$S$};
    \filldraw (S.180) circle (.5mm);
\end{tikzpicture}
= i \cdot
\begin{tikzpicture}[baseline]
    \clip (-.6,-1) rectangle (.6,1);

    \draw(-.32,0) -- ++(90:15mm);
    \draw (.32,0) -- ++(90:15mm);
    \draw (-.07,.7) node {...};
    \draw (.15,0) -- ++(90:15mm);
    
    \node (S) at (0,0) [circle, draw, fill=white] {$S$};
    \filldraw (S.180) circle (.5mm);
\end{tikzpicture}
$

%% file: diagrams/tikz/capS.tex
$
\begin{tikzpicture}[baseline]
 \node (S) [circle,draw] {$S$};
 \filldraw (S.157) circle (.5mm);
 \draw (S.135) .. controls +(135:3mm) and +(90:3mm) .. (S.90);
 \draw (S.45) -- +(45:3mm);
 \draw (S.0) -- +(0:3mm);
 \draw (S.-45) node [below] {$\cdot$};
  \draw (S.-90) node [below] {$\cdot$};
 \draw (S.-135) node [below] {$\cdot$};
 \draw (S.180) -- +(180:3mm);
\end{tikzpicture}
= 0
$

%% file: diagrams/tikz/twoS.tex
$
\begin{tikzpicture}[baseline]
 \node at (0,.7) (S) [circle,draw] {$S$};
     \filldraw (S.180) circle (.5mm);

 \node at (0,-.7) (S') [circle,draw] {$S$};
     \filldraw (S'.180) circle (.5mm);

 \draw (S.45) -- +(45:4mm);
 \draw (S.90) node [above] {$\ldots$};
 \draw (S.135) -- +(135:4mm);
 \draw (S'.225) -- +(-135:4mm);
 \draw (S'.270) node [below] {$\ldots$};
 \draw (S'.-45) -- +(-45:4mm);
\end{tikzpicture}
=\qi{2n-1} \cdot
\begin{tikzpicture}[baseline]
\node (0,0) (JW) [rectangle,draw] {$\JW{4n-4}$};
\draw (JW.180);

\draw (JW.35) -- +(90:11mm);
\draw (JW.90) node [above] {$\cdots$};
\draw (JW.145) -- +(90:11mm);

\draw (JW.-35) -- +(-90:11mm);
\draw (JW.-90) node [below] {$\cdots$};
\draw (JW.-145) -- +(-90:11mm);
\end{tikzpicture}
$

%% file: text/background.tex
In this section we remind the reader what a planar algebra is, and recall a few facts about the simplest planar algebra, Temperley-Lieb.

\subsection{What is a planar algebra?}\label{subsec:pa}
A planar algebra is a gadget specifying how to combine elements in planar
ways, rather as a ``linear" algebra is a gadget in which one can combine elements, given a linear ordering.  For  example,
$$ f \circ g \circ h \text{ \quad versus \quad }\mathfig{0.2}{tangles/cubic-fgh}.$$
Planar algebras were introduced in \cite{math.QA/9909027} to study subfactors,
and have since found more general use.

In the simplest version, a planar algebra $\mathcal{P}$ associates a
vector space $\mathcal{P}_k$ to each natural number $k$ (thought of as a
disc in the plane with $k$ points on its boundary) and associates a linear map
$\mathcal{P}(T)  : \mathcal{P}_{k_1} \tensor \mathcal{P}_{k_2} \tensor
\cdots \tensor \mathcal{P}_{k_r} \to \mathcal{P}_{k_0}$ to each
`spaghetti and meatballs' or `planar tangle' diagram $T$ with internal
discs with $k_1, k_2, \ldots, k_r$ points and $k_0$ points on the
external disc.  For example,
$$\mathfig{0.2}{tangles/cubic}$$ gives a map from $V_7 \tensor V_5 \tensor V_5 \rightarrow V_7$. Such maps (the `planar operations') must satisfy certain
properties: radial spaghetti induces the identity map, and composition of the
maps $\mathcal{P}(T)$ is compatible with the obvious composition of
spaghetti and meatballs diagrams by gluing some inside another.  When we glue, we match up base points;
each disk's base point is specified by a bullet.

The reason for these bullets in the definition is that they allow us to keep track of pictures which are not rotationally invariant.  For example, in Definition \ref{def:pa} we have used the marked points to indicate the way the generator $S$ behaves under rotation. Nevertheless we use the following conventions to avoid always drawing a bullet.
Instead of using marked points we will often instead use a ``rectangular" picture in which some of the strings go up, some go down, but none come out of the sides of generators.  This leaves a gap on the left side of every picture and the convention is that the marked points always lie in this gap.  Further, if we neglect to draw a bounding disk, the reader should imagine a rectangle around the picture (and therefore put the marked point on the left).  For example, the following pictures represent the same element of a planar algebra: $S$ rotated by one `click'.

\begin{center}
\begin{tikzpicture}[baseline]
    \node (S) [circle,draw] {$S$};
    \draw (0,0) circle (1cm);
    \filldraw (S.157) circle (.5mm);
    \draw (S.180)-- (180:10mm);
    \draw (S.135)-- (135:10mm);
    \draw (S.90)-- (90:10mm);
    \draw (S.45)-- (45:10mm);
    \draw (S.0)-- (0:10mm);
    \draw (S.-45)-- (-45:10mm);
    \filldraw (112:1cm) circle (.5mm);
    \draw (S.-90)-- (-90:10mm);
    \draw (S.-135)-- (-135:10mm);
\end{tikzpicture}
\quad , \quad
\begin{tikzpicture}[baseline]
    \clip (-1.1,-1.4) rectangle (1.1,1.4);

    \draw[rounded corners=2mm] (-.35,0) -- ++(90:.7cm) -- ++(180:.5cm) -- ++ (-90:1.7cm) -- ++(0:1.6cm) -- ++(90:3cm);
    \draw (.35,0) -- ++(90:15mm);
    \draw (-.25,0) -- ++(90:15mm);
    \draw (-.15,0) -- ++(90:15mm);
    \draw (-.05,0) -- ++(90:15mm);
    \draw (.05,0) -- ++(90:15mm);
    \draw (.15,0) -- ++(90:15mm);
    \draw (.25,0) -- ++(90:15mm);
    
    \node (S) at (0,0) [rectangle, draw, fill=white] {$\;\; S \;\;$};
\end{tikzpicture}
\end{center}

There are some special planar tangles which
induce operations with familiar names. First, each even vector space $\mathcal{P}_{2k}$ becomes an associative algebra using the `multiplication' tangle:
$$\mathfig{0.2}{tangles/multiplication}$$
Second, there is an involution $\overline{\phantom{X}} : \mathcal{P}_{2k} \to \mathcal{P}_{2k}$ given by the `dualising' tangle:
$$\mathfig{0.2}{tangles/dualising}$$
Third, for each $k$ there is a trace $\text{tr} : \mathcal{P}_{2k} \rightarrow \mathcal{P}_0$:
$$\mathfig{0.2}{tangles/trace}$$
If $\mathcal{P}_0$ is one-dimensional, this map really is a trace, and we can use it (along with multiplication) to build a bilinear form on $P_{2k}$ in the usual way.

A subfactor planar algebra is the best kind of planar algebra; it has additional properties which make it a nice place to work.  First and foremost, $\mathcal{P}_0$ must be one-dimensional.  In particular, a closed circle is equal to a multiple of the empty diagram, and the square of this multiple is called the \emph{index} of the planar algebra.  Note that this implies that the zero-ary planar operations, namely the `vegetarian' diagrams without any meatballs, induce the Temperley-Lieb diagrams (see below, \S \ref{ssTL}) as elements of the subfactor planar algebra.
There is thus a map $\TL_{2k} \to \mathcal{P}_{2k}$, although it need be neither surjective nor injective.

Second, subfactor planar algebras have the property that only spaces for discs with an even number of boundary points are nonzero. Third, subfactor planar algebras must be spherical, that is, for each element $T \in \mathcal{P}_2$, we have an identity in $\mathcal{P}_0$:
\begin{equation*}
\mathfig{0.1}{spherical/1} = \mathfig{0.1}{spherical/3}.
\end{equation*}
Fourth, there must be an anti-linear \emph{adjoint} operation $*:\mathcal{P}_k \to \mathcal{P}_k$ such that the sesquilinear form given by $\left< x , y \right> = \tr{y^*x}$ is positive definite.  Further, $*$ on $\mathcal{P}$ should be compatible with the horizontal reflection operation $*$ on planar tangles.  In particular, this means that the adjoint operation on Temperley-Lieb is reflection in a horizontal line.

Finally note that we use ``star'' to indicate the \emph{adjoint}, and ``bar'' to indicate the \emph{dual}. We apologise
to confused readers for this notation.

One useful way to generalize the definition of a planar algebra is by introducing a `spaghetti label set' and a `region label set,'
and perhaps insist that only certain labels can appear next to each other. When talking about subfactor planar algebras, only two simple cases of this
are required: a `standard' subfactor planar algebra has just two region labels, shaded and unshaded, which must alternate across spaghetti, while an `unshaded' subfactor planar algebra has no interesting labels at all.

From this point onwards, we'll be using the unshaded variety of planar algebra, essentially just for simplicity of exposition. The reader can easily
reconstruct the \emph{shaded} version of everything we say; checkerboard shade the regions, ensuring that the marked point of an $S$ box is always in an unshaded region.
This necessitates replacing relation \ref{rotateS} in Definition \ref{def:pa}, so that instead the ``2 click'' rotation of the $S$ box is $-1$ times the original unrotated box.
The one point at which reintroducing the shading becomes subtle is when we discuss braidings in \S \ref{sec:braiding}.

\subsection{The Temperley-Lieb (planar) algebra} \label{ssTL}
We work over the field $\Complex(q)$ of rational functions in a formal variable $q$.
It is often notationally convenient to use quantum numbers.

\begin{defn}
The $n^{th}$ quantum number $\qi{n}$ is defined as $$\frac{q^n-q^{-n}}{q-q^{-1}} = q^{n-1} + q^{n-3} + \cdots + q^{-n+1}.$$
 \end{defn}

Now let's recall some facts about the Temperley-Lieb algebra.

\begin{defn}
A Temperley-Lieb picture is a non-crossing matching of $2n$ points around the boundary of a disc, with a chosen first point.
\end{defn}

In practice, Temperley-Lieb pictures are often drawn with the points on two lines, and the chosen first point is the one on the top left.

\begin{defn}
The vector space $\TL_{2n}$ has as a basis the Temperley-Lieb pictures for $2n$ points.  These assemble into a planar
algebra by gluing diagrams into planar tangles, and removing each
closed circle formed in exchange for a coefficient of $\qi{2} = q + q^{-1}$.
\end{defn}

Temperley-Lieb is a subfactor planar algebra (with the adjoint operation being horizontal reflection) except that the sesquilinear form need not be positive definite (see \S \ref{sec:special-TL}).  Some important elements of the Temperley-Lieb algebra are

\begin{itemize}
\item The identity (so-called because it's the identity for the multiplication given by vertical stacking):
$$\id =\input{diagrams/tikz/Identity.tex};$$
\item The Jones projections in $\TL_{2n}$:
$$e_i=\input{diagrams/tikz/JonesProjection.tex}, \, i \in \{ 1, \ldots, n-1\};$$
\item
The Jones-Wenzl projection \cite{MR873400} $\JW{n}$ in $\TL_{2n}$:  The unique projection
with the property
$$\JW{n} e_i = e_i \JW{n} =0, \, \text{for each $i \in \{1, \ldots, n-1 \}$};$$
\item The crossing in $\TL_4$:
\begin{equation*}
\input{diagrams/tikz/ResolvedCrossing.tex}.
\end{equation*}
\end{itemize}

Recall that the crossing satisfies Reidemeister relations $2$ and $3$, but not Reidemeister $1$.  Instead the positive twist factor is $i q^{3/2}$.

\begin{fact}
Here are some useful identities involving the Jones-Wenzl projections:
\begin{enumerate}
\item
{ \bf Wenzl's relation:}
\begin{equation} \label{WenzlsRelation} \input{diagrams/tikz/WenzlRelation.tex} \end{equation}

and $\tr{\JW{m}} = \qi{m+1}$.

\item
{\bf Partial trace relation:} \begin{equation}\label{eq:PartialTraceRelation}\input{diagrams/tikz/PartialTraceRelation.tex}\end{equation}

\end{enumerate}
\end{fact}

\subsection{Temperley-Lieb when $f^{(4n-3)} = 0$}
\label{sec:special-TL}%
At any `special value' $q=e^{\frac{i \pi}{k+2}}$ (equivalently $\delta = q+q^{-1} = 2 \cos(\frac{\pi}{k+2})$), the
Temperley-Lieb planar algebra is degenerate, with radical generated by
the Jones-Wenzl projection $f^{(k+1)}$. We therefore pass to a
quotient, by imposing the relation $f^{(k+1)} = 0$. In the physics literature $k$ would be called the level.  We're interested in the case $k=4n-4$, so
$q=e^{i \pi/(4n-2)}$ and $\delta=2 \cos{ \frac{\pi}{4n-2}}$.  For this value of $q$, $\qi{m}=\qi{4n-2-m}$.

We record several facts about this quotient of Temperley-Lieb which we'll need later.
(In the following diagrams, we're just drawing $3$ or $4$ parallel strands where we really mean $4n-5$ or $4n-4$ respectively; make sure you read the labels of the boxes.)

\begin{lem} Strands cabled by $f^{(4n-4)}$ can be reconnected.
\label{lem:2f}
$$\mathfig{0.3}{consistency/2f} = \rotatemathfig{0.3}{90}{consistency/2f}$$
\end{lem}
\begin{rem}
Any relation in Temperley-Lieb also holds if superimposed on top of, or
behind, another Temperley-Lieb diagram; this is just the statement that
Temperley-Lieb is braided. We'll need to use all these variations of the
identity in the above lemma later.
\end{rem}

\begin{lem}
\label{lem:f-twist}
\begin{align*}
\rotatemathfig{0.15}{-90}{consistency/f-twist-1} & =
\rotatemathfig{0.15}{-90}{consistency/f-twist-2}
     =
     q^{2n-1} \rotatemathfig{0.12}{-90}{consistency/f-no-twist}
     =
     i \rotatemathfig{0.12}{-90}{consistency/f-no-twist}
\end{align*}
(the twisted strand here indicates just a single strand, while the 3 parallel strands actually represent 4n-5 strands)
and as an easy consequence
$$\rotatemathfig{0.25}{-90}{consistency/f-twist-full} =
\rotatemathfig{0.12}{-90}{consistency/f-no-twist}.$$
\end{lem}

The first two equalities hold in Temperley-Lieb at any value of $q$. The
third equality simply specialises to the relevant value. Note that the
crossings in the above lemma are all undercrossings for the single strand. Changing each
of these to an overcrossing for that strand, we have the same identities, with $q$ replaced by
$q^{-1}$, and $i$ replaced by $-i$.

\begin{lem} Overcrossings, undercrossings and the $2$-string identity cabled by $f^{(4n-4)}$ are all the same.
\label{lem:f-crossings}
$$\mathfig{0.3}{consistency/f-crossings-1} = \mathfig{0.3}{consistency/f-crossings-2} =
\mathfig{0.3}{consistency/2f}$$
\end{lem}

%% file: diagrams/tikz/Identity.tex
\begin{tikzpicture}[baseline]
    \draw (.2,-.5)--(.2,.5);
	\draw (.4,-.5) -- (.4,.5);
	\node at (.8,0) {$\cdots$};
    \draw (1.1,-.5)--(1.1,.5);
\end{tikzpicture}

%% file: diagrams/tikz/JonesProjection.tex
\qi{2}^{-1} \begin{tikzpicture}[baseline]
    \draw (.2,-.5)--(.2,.5);
    \node at (.6,0) {$\cdots$};
    \draw (1,-.5)--(1,.5);
    \draw (1.2,-.5) arc (180:0:1mm);
    \draw (1.2,.5) arc (-180:0:1mm);
    \draw (1.6,-.5) -- (1.6,.5);
    \node at (2,0) {$\cdots$};
    \draw (2.3,-.5) -- (2.3,.5);
\end{tikzpicture}

%% file: diagrams/tikz/ResolvedCrossing.tex
\begin{tikzpicture}[baseline]
\node (x) at (0,0){};
    \draw (x.45)-- (.5,.5);
    \draw (x.135) -- (-.5,.5);
    \draw (x.315) -- (.5,-.5);
    \draw (x.45) -- (-.5,-.5);
\end{tikzpicture}
:= i q^{\frac 1 2}
\begin{tikzpicture}[baseline]
    \draw (1.5,.5) .. controls (2,0) .. (1.5,-.5);
    \draw (2.5,.5) .. controls (2,0) .. (2.5,-.5);
\end{tikzpicture}
-i q^{-\frac 1 2}
\begin{tikzpicture}[baseline]
    \draw (3.5,.5) .. controls (4,0) .. (4.5,.5);
    \draw (3.5,-.5) .. controls (4,0) .. (4.5,-.5);
\end{tikzpicture}

%% file: diagrams/tikz/WenzlRelation.tex
\qi{m+1} \cdot
\begin{tikzpicture}[baseline]
    \clip (-.9,-1.3) rectangle (1,1.3);
    \node at (0,0) (box) [rectangle,draw] {$\hspace{.1cm} \JW{m + 1} \hspace{.1cm} $};

    \draw (box.-50)-- ++(-90:1cm);
    \draw (box.-90) -- ++(-90:1cm);
    \draw (box.-150)--++(-90:1cm);
    \node at (-.25,-.6) {...};
    \draw (box.-30) -- ++(-90:1cm);

    \draw (box.50)-- ++(90:1cm);
    \draw (box.90) -- ++(90:1cm);
    \draw (box.150)--++(90:1cm);
    \node at (-.25,.6) {...};
    \draw (box.30) -- ++(90:1cm);
\end{tikzpicture}
=
\qi{m+1} \cdot
\begin{tikzpicture}[baseline]
    \clip (-.8,-1.3) rectangle (1,1.3);
    \node at (0,0) (box) [rectangle,draw] {\; $\JW{m}$ \;};

    \draw (box.-145)--++(-90:1cm);
    \node at (-.1,-.6) {...};
    \draw (box.-35) -- ++(-90:1cm);
    \draw (box.-60) -- ++(-90:1cm);

    \draw (box.35) -- ++(90:1cm);
    \draw (box.60) -- ++(90:1cm);
    \node at (-.1,.6) {...};
    \draw (box.145) -- ++(90:1cm);

    \draw (.9,1.5)--(.9,-1.5);
\end{tikzpicture}
-
\qi{m} \cdot
\begin{tikzpicture}[baseline]
    \clip (-.8,-1.3) rectangle (1,1.3);
    \node at (0,.7) (top) [rectangle,draw] {$\;\; \JW{m} \; \;$};
    \node at (0,-.7) (bottom) [rectangle,draw]{$\; \;\JW{m} \; \;$};

    \draw (top.-60)--(bottom.60);
    \draw (top.-145)--(bottom.145);
    \node at (-.1,0) {...};
    \draw (top.-35) arc (-180:0:.2cm) -- ++(90:1cm);
    \draw (bottom.35) arc (180:0:.2cm) -- ++(-90:1cm);

    \draw (top.35) -- ++(90:1cm);
    \draw (top.60) -- ++(90:1cm);
    \node at (-.1,1.2) {...};
    \draw (top.145) -- ++(90:1cm);

    \draw (bottom.-35) -- ++(-90:1cm);
    \draw (bottom.-60) -- ++(-90:1cm);
    \node at (-.1,-1.2) {...};
    \draw (bottom.-145) -- ++(-90:1cm);
\end{tikzpicture}

%% file: diagrams/tikz/PartialTraceRelation.tex
\qi{m} \cdot
\begin{tikzpicture}[baseline]
    \node (0,0) (JW) [rectangle,draw] {$\; \; \JW{m} \; \;$};

    \draw [rounded corners=2mm] (JW.35) -- ++(90:4mm) -- ++(4mm,0mm) -- ++(0cm,-1.5cm) -- ++(-4mm,  0mm)-- (JW.-35);
    \draw (JW.60) -- +(90:5mm);
    \draw (JW. 110) node [above] {$\cdots$};
    \draw (JW.145) -- +(90:5mm);

    \draw (JW.-60) -- +(-90:5mm);
    \draw (JW.-110) node [below] {$\cdots$};
    \draw (JW.-145) -- +(-90:5mm);
\end{tikzpicture}
= \qi{m+1} \cdot
\begin{tikzpicture}[baseline]
    \node (0,0) (JW) [rectangle,draw] {$\; \JW{m-1} \;$};

    \draw (JW.35) -- +(90:5mm);
    \draw (JW. 90) node [above] {$\cdots$};
    \draw (JW.145) -- +(90:5mm);

    \draw (JW.-35) -- +(-90:5mm);
    \draw (JW.-90) node [below] {$\cdots$};
    \draw (JW.-145) -- +(-90:5mm);
\end{tikzpicture}

%% file: text/presentation.tex
\subsection{First consequences of the relations}
Recall from the introduction that we are considering the planar algebra $\pa$ generated by a single box $S$ with $4n-4$ strands, with $q=\exp(\frac{\pi i}{4n-2})$, modulo the following relations.
\begin{enumerate}

\item  A closed circle is equal to $[2]_q = (q+q^{-1}) = 2 \cos(\frac{\pi}{4n-2})$ times the empty diagram.

\item\input{diagrams/tikz/rotateS.tex}

\item\input{diagrams/tikz/capS.tex}

\item\input{diagrams/tikz/twoS.tex}

\end{enumerate}

\begin{rem}
Relation \eqref{delta} fixes the index $\qi{2}^2$ of the planar algebra as a `special value' as in \S \ref{sec:special-TL} of the form $\qi{2} = 2 \cos(\frac{\pi}{k+2})$.
Note that usually at special values, one imposes a further relation, that the corresponding Jones-Wenzl idempotent $f^{(k)}$ is zero, in order that the planar algebra be positive definite.
As it turns out, we \emph{don't} need to impose this relation by hand; it will follow, in Theorem \ref{easyconsequences}, from the other relations.

According to the philosophy of \cite{MR1929335, quadratic} any planar algebra is generated by boxes which satisfy ``annular relations"  like \eqref{rotateS} and \eqref{capS}, while particularly nice planar algebras require in addition only ``quadratic relations" which involve two boxes.  Our quadratic relation \eqref{twoS}, in which the two $S$ boxes are not connected, is unusually strong and makes many of our subsequent arguments possible. Notice that this relation also implies relations with a pair of $S$ boxes connected by an arbitrary number of strands.
\end{rem}

We record for future use some easy consequences of the relations of Definition \ref{def:pa}.

\begin{thm}\label{easyconsequences}
The following relations hold in $\pa$.
\begin{enumerate}
\item
\input{diagrams/tikz/StimesS.tex}

(Here $2n-2$ strands connect the two $S$ boxes on the left hand side.)
\item
\input{diagrams/tikz/otherScaps.tex}

\item\label{JWzero}
\input{diagrams/tikz/JWzerostatement.tex}

\item\label{enoughcappings}
 For $T, T'  \in \pa_{4n-3}$, if

$$
\input{diagrams/tikz/enoughcappings.tex}
$$

then $T=T'$.  More generally, if $T,T' \in \pa_m$ for $m \geq {4n-3}$, and $4n-4$ consecutive cappings of $T$ and $T'$ are equal, then $T=T'$.

\end{enumerate}
\end{thm}

\begin{proof}
\begin{enumerate}

\item This follows from taking a partial trace (that is, connecting top right strings to bottom right strings) of the diagrams of the two-$S$ relation (\ref{twoS}), and applying the partial trace relation from Equation \eqref{eq:PartialTraceRelation}.

\item This is a straightforward application of the rotation relation (\ref{rotateS}) and the capping relation (\ref{capS}).

\item Using Wenzl's relation (Equation \eqref{WenzlsRelation}) we calculate

\begin{align*}
\input{diagrams/tikz/JWzeroWenzl.tex}
\intertext{then replace $\qi{4n-3}$ and $\qi{4n-4}$ by $1$ and $\qi{2}$, and apply the two $S$ relation thrice, obtaining}
\input{diagrams/tikz/JWzeroputinS.tex}
\intertext{We can then use the two-$S$ relation on the middle two $S$ boxes of the second picture, and apply the partial trace relation (Equation \eqref{eq:PartialTraceRelation}) to the resulting $\JW{4n-4}$. We thus see}
\qi{4n-3}
\begin{tikzpicture}[baseline]
    \clip (-1,-1.3) rectangle (1,1.3);
    \node at (0,0) (box) [rectangle,draw] {$\hspace{.1cm} \JW{4n -3} \hspace{.1cm} $};
    \draw (box.-50)-- ++(-90:1cm);
    \draw (box.-90) -- ++(-90:1cm);
    \draw (box.-150)--++(-90:1cm);
    \node at (-.25,-.6) {...};
    \draw (box.-30) -- ++(-90:1cm);
    \draw (box.50)-- ++(90:1cm);
    \draw (box.90) -- ++(90:1cm);
    \draw (box.150)--++(90:1cm);
    \node at (-.25,.6) {...};
    \draw (box.30) -- ++(90:1cm);
\end{tikzpicture} \input{diagrams/tikz/JWzerotakeoutS.tex}
\input{diagrams/tikz/JWzerosimplify.tex}
\end{align*}

\item Thanks to Stephen Bigelow for pointing out this fact.

On the one hand, $\JW{4n-3}$ is a weighted sum of Temperley-Lieb pictures, with the weight of $\id$ being $1$:   $$\JW{4n-3} = \id + \sum_{P \in \TL_{4n-3}, P \neq \id} \alpha_P \cdot P;$$

On the other hand, $\JW{4n-3}=0$.  Therefore $$\id  = \id - \JW{4n-3} = \sum_{P \in \TL_{4n-3},  P \neq \id} -\alpha_P \cdot P.$$
If $P \in \TL_{4n-3}$ and $P \neq \id$, then $P$ has a cap somewhere along the boundary, so it follows from our hypotheses that $P T = P T'$, and therefore $$T = \left ( \sum_{P \in \TL_{4n-3},  P \neq \id} -\alpha_P \cdot P\right) T =  \left( \sum_{P \in \TL_{4n-3},  P \neq \id} -\alpha_P \cdot P \right) T' = T' .$$
\end{enumerate}

\end{proof}

%% file: diagrams/tikz/StimesS.tex
$
\begin{tikzpicture}[baseline]
	\clip (-.8,-1.3) rectangle (.8,1.3);
	 \node at (0,.6) (S) [circle,draw] {$ S$};
	 \filldraw (S.180) circle (.5mm);
	\node at (0,-.6) (S') [circle,draw] {$S$};
	 \filldraw (S'.180) circle (.5mm);

     \draw (S.203)..controls +(203:2mm) and +(157:2mm) .. (S'.157);
     \draw (S'.90) node[above] {$\ldots$};
     \draw (S.337)..controls +(337:2mm) and +(23:2mm) ..
     (S'.23);
     \draw (S.45) -- +(45:1cm);
     \draw (S.90) node [above] {$\ldots$};
     \draw (S.135) -- +(135:1cm);
     \draw (S'.225) -- +(-135:1cm);
     \draw (S'.270) node [below] {$\ldots$};
     \draw (S'.-45) -- +(-45:1cm);
\end{tikzpicture}
=
\begin{tikzpicture}[baseline]
	\clip (-.8,-1.3) rectangle (.8,1.3);
    \node (0,0) (JW) [rectangle,draw] {$\JW{2n-2}$};

    \draw (JW.35) -- +(90:11mm);
    \draw (JW.90) node [above] {$\cdots$};
    \draw (JW.145) -- +(90:11mm);

    \draw (JW.-35) -- +(-90:11mm);
    \draw (JW.-90) node [below] {$\cdots$};
    \draw (JW.-145) -- +(-90:11mm);
\end{tikzpicture}
$

%% file: diagrams/tikz/otherScaps.tex
$
\begin{tikzpicture}[baseline]
 \node (S) [circle,draw] {$ S$};
  \filldraw (S.157) circle (.5mm);
     \draw (S.135) -- +(135:3mm);
     \draw (S.90) -- +(90:3mm);
     \draw (S.45) .. controls +(45:3mm) and +(0:3mm) .. (S.0);
     \draw (S.-45) node [below] {$\cdot$};
      \draw (S.-90) node [below] {$\cdot$};
     \draw (S.-135) node [below] {$\cdot$};
     \draw (S.180) -- +(180:3mm);
\end{tikzpicture}
=
\begin{tikzpicture}[baseline]
\node[xshift=2cm] (S') [circle,draw] {$ S$};
 \filldraw (S'.157) circle (.5mm);
     \draw (S'.135) -- +(135:3mm);
     \draw (S'.90) -- +(90:3mm);
     \draw (S'.45) -- +(45:3mm);
     \draw (S'.0) .. controls +(0:3mm) and +(-45:3mm) .. (S'.-45);
      \draw (S'.-90) node [below] {$\cdot$};
     \draw (S'.-120) node [below] {$\cdot$};
     \draw (S'.-150) node [below] {$\cdot$};
     \draw (S'.180) -- +(180:3mm);
\end{tikzpicture}
= \cdots =0
$

%% file: diagrams/tikz/JWzerostatement.tex
$\begin{tikzpicture}[baseline]
    \node (0,0) (JW) [rectangle,draw] {$\JW{4n-3}$};

    \draw (JW.35) -- +(90:5mm);
    \draw (JW.90) node [above] {$\cdots$};
    \draw (JW.145) -- +(90:5mm);

    \draw (JW.-35) -- +(-90:5mm);
    \draw (JW.-90) node [below] {$\cdots$};
    \draw (JW.-145) -- +(-90:5mm);
\end{tikzpicture}
=0$

%% file: diagrams/tikz/enoughcappings.tex
\begin{tikzpicture}[baseline]
	\draw (-.6,0) -- (-.6,.3) arc (180:0:1mm) -- (-.4,0);
	\draw (-.2,0) -- (-.2,.5);
	\draw (.4,0) -- (.4,.5);
	\draw (.6,0) -- (.6,.5);	
	\node at (.1,.4) {...};
	\node at (0,0) (T) [rectangle, draw, fill=white] {$\hspace{4mm} T \hspace{4mm}$};
\end{tikzpicture}
=
\begin{tikzpicture}[baseline]
	\draw (-.6,0) -- (-.6,.3) arc (180:0:1mm) -- (-.4,0);
	\draw (-.2,0) -- (-.2,.5);
	\draw (.4,0) -- (.4,.5);
	\draw (.6,0) -- (.6,.5);	
	\node at (.1,.4) {...};
	\node at (0,0) (T) [rectangle, draw, fill=white] {$\hspace{4mm} T' \hspace{4mm}$};
\end{tikzpicture}
\, , \,
\begin{tikzpicture}[baseline]
	\draw (-.6,0) -- (-.6,.5);
	\draw (-.4,0) -- (-.4,.3) arc (180:0:1mm) -- (-.2,0);
	\draw (.4,0) -- (.4,.5);
	\draw (.6,0) -- (.6,.5);	
	\node at (.1,.4) {...};
	\node at (0,0) (T) [rectangle, draw, fill=white] {$\hspace{4mm} T \hspace{4mm}$};
\end{tikzpicture}
=\begin{tikzpicture}[baseline]
	\draw (-.6,0) -- (-.6,.5);
	\draw (-.4,0) -- (-.4,.3) arc (180:0:1mm) -- (-.2,0);
	\draw (.4,0) -- (.4,.5);
	\draw (.6,0) -- (.6,.5);	
	\node at (.1,.4) {...};
	\node at (0,0) (T') [rectangle, draw, fill=white] {$\hspace{4mm} T' \hspace{4mm}$};
\end{tikzpicture}
\, , \ldots , \,
\begin{tikzpicture}[baseline]
	\draw (-.6,0) -- (-.6,.5);
	\draw (-.4,0) -- (-.4,.5);
	\draw (-.2,0) -- (-.2,.5);
	\draw (.4,0) -- (.4,.3)  arc (180:0:1mm) -- (.6,0);
	\node at (.1,.4) {...};
	\node at (0,0) (T) [rectangle, draw, fill=white] {$\hspace{4mm} T \hspace{4mm}$};
\end{tikzpicture}
=\begin{tikzpicture}[baseline]
	\draw (-.6,0) -- (-.6,.5);
	\draw (-.4,0) -- (-.4,.5);
	\draw (-.2,0) -- (-.2,.5);
	\draw (.4,0) -- (.4,.3)  arc (180:0:1mm) -- (.6,0);
	\node at (.1,.4) {...};
	\node at (0,0) (T') [rectangle, draw, fill=white] {$\hspace{4mm} T' \hspace{4mm}$};
\end{tikzpicture}

%% file: diagrams/tikz/JWzeroWenzl.tex
\qi{4n-3}
\begin{tikzpicture}[baseline]
    \clip (-1,-1.3) rectangle (1,1.3);
    \node at (0,0) (box) [rectangle,draw] {$\hspace{.1cm} \JW{4n -3} \hspace{.1cm} $};

    \draw (box.-50)-- ++(-90:1cm);
    \draw (box.-90) -- ++(-90:1cm);
    \draw (box.-150)--++(-90:1cm);
    \node at (-.25,-.6) {...};
    \draw (box.-30) -- ++(-90:1cm);

    \draw (box.50)-- ++(90:1cm);
    \draw (box.90) -- ++(90:1cm);
    \draw (box.150)--++(90:1cm);
    \node at (-.25,.6) {...};
    \draw (box.30) -- ++(90:1cm);
\end{tikzpicture}
& = \qi{4n-3}
\begin{tikzpicture}[baseline]
    \clip (-1,-1.3) rectangle (1,1.3);
    \node at (0,0) (box) [rectangle,draw] {$\JW{4n-4}$};

    \draw (box.-145)--++(-90:1cm);
    \node at (-.1,-.6) {...};
    \draw (box.-35) -- ++(-90:1cm);
    \draw (box.-60) -- ++(-90:1cm);

    \draw (box.35) -- ++(90:1cm);
    \draw (box.60) -- ++(90:1cm);
    \node at (-.1,.6) {...};
    \draw (box.145) -- ++(90:1cm);

    \draw (.9,1.5)--(.9,-1.5);
\end{tikzpicture}
-
\qi{4n-4}
\begin{tikzpicture}[baseline]
    \clip (-1,-1.3) rectangle (1,1.3);
    \node at (0,.7) (top) [rectangle,draw] {$\JW{4n-4}$};
    \node at (0,-.7) (bottom) [rectangle,draw]{$\JW{4n-4}$};

    \draw (top.-60)--(bottom.60);
    \draw (top.-145)--(bottom.145);
    \node at (-.1,0) {...};
    \draw (top.-35) arc (-180:0:.2cm) -- ++(90:1cm);
    \draw (bottom.35) arc (180:0:.2cm) -- ++(-90:1cm);

    \draw (top.35) -- ++(90:1cm);
    \draw (top.60) -- ++(90:1cm);
    \node at (-.1,1.2) {...};
    \draw (top.145) -- ++(90:1cm);

    \draw (bottom.-35) -- ++(-90:1cm);
    \draw (bottom.-60) -- ++(-90:1cm);
    \node at (-.1,-1.2) {...};
    \draw (bottom.-145) -- ++(-90:1cm);
\end{tikzpicture}, \\

%% file: diagrams/tikz/JWzeroputinS.tex
& = \frac{1}{\qi{2n-1}}
\begin{tikzpicture}[baseline]
    \clip (-1,-1.3) rectangle (1,1.3);

    \node at (0,.6) (top) [ellipse,draw] {$ \hspace{.15 cm} S \hspace{.15cm} $};
    \node at (0,-.6) (bottom) [ellipse,draw]{$ \hspace{.15 cm} S \hspace{.15cm} $};
    \filldraw (top.180) circle (.5mm);
    \filldraw (bottom.180) circle (.5mm);    

    \node at (top.0) (invisible) [right] {};
    \draw (invisible.90)--++(90:2cm);
    \draw (invisible.90) -- ++ (-90:2cm);

    \draw (top.35) -- ++(90:1cm);
    \draw (top.60) -- ++(90:1cm);
    \node at (-.1,1.2) {...};
    \draw (top.145) -- ++(90:1cm);

    \draw (bottom.-35) -- ++(-90:1cm);
    \draw (bottom.-60) -- ++(-90:1cm);
    \node at (-.1,-1.2) {...};
    \draw (bottom.-145) -- ++(-90:1cm);
\end{tikzpicture}
-
\frac{\qi{2}}{\qi{2n-2}^2}
\begin{tikzpicture}[baseline]
    \clip (-1,-2.2) rectangle (1,2.2);

    \node at (0,.6) (top) [ellipse,draw] {$ \hspace{.15 cm} S \hspace{.15cm} $};
    \node at (0,-.6) (bottom) [ellipse,draw]{$ \hspace{.15 cm} S \hspace{.15cm} $};
    \node at (0,-1.6) (bottom2) [ellipse,draw]{$ \hspace{.15 cm} S \hspace{.15cm} $};
    \node at (0,1.6) (top2) [ellipse,draw]{$ \hspace{.15 cm} S \hspace{.15cm} $};
    \filldraw (top.180) circle (.5mm);
    \filldraw (bottom.180) circle (.5mm);        
    \filldraw (top2.180) circle (.5mm);
    \filldraw (bottom2.180) circle (.5mm);    

    \draw (top.-60)--(bottom.60);
    \draw (top.-145)--(bottom.145);
    \node at (-.1,0) {...};
    \draw (top.-35) arc (-180:0:.2cm) -- ++(90:2cm);
    \draw (bottom.35) arc (180:0:.2cm) -- ++(-90:2cm);

    \draw (top2.35) -- ++(90:1cm);
    \draw (top2.60) -- ++(90:1cm);
    \node at (-.1,2.1) {...};
    \draw (top2.145) -- ++(90:1cm);

    \draw (bottom2.-35) -- ++(-90:1cm);
    \draw (bottom2.-60) -- ++(-90:1cm);
    \node at (-.1,-2.1) {...};
    \draw (bottom2.-145) -- ++(-90:1cm);
\end{tikzpicture}, \\

%% file: diagrams/tikz/JWzerotakeoutS.tex
&= \frac{1}{\qi{2n-1}}
\begin{tikzpicture}[baseline]
    \clip (-1,-1.3) rectangle (1,1.3);

    \node at (0,.6) (top) [ellipse,draw] {$ \hspace{.15 cm} S \hspace{.15cm} $};
    \node at (0,-.6) (bottom) [ellipse,draw]{$ \hspace{.15 cm} S \hspace{.15cm} $};
    \filldraw (top.180) circle (.5mm);
    \filldraw (bottom.180) circle (.5mm);    

    \node at (top.0) (invisible) [right] {};
    \draw (invisible.90)--++(90:2cm);
    \draw (invisible.90) -- ++ (-90:2cm);

    \draw (top.35) -- ++(90:1cm);
    \draw (top.60) -- ++(90:1cm);
    \node at (-.1,1.2) {...};
    \draw (top.145) -- ++(90:1cm);

    \draw (bottom.-35) -- ++(-90:1cm);
    \draw (bottom.-60) -- ++(-90:1cm);
    \node at (-.1,-1.2) {...};
    \draw (bottom.-145) -- ++(-90:1cm);
\end{tikzpicture}
-
\frac{\qi{2}}{\qi{2n-1}}
\begin{tikzpicture}[baseline]
    \clip (-1,-2.2) rectangle (2.5,2.2);

    \node at (.5,0) (JW) [rectangle, draw] {$ \JW{4n-4}$};
    \node at (0,-1.6) (bottom2) [ellipse,draw]{$ \hspace{.15 cm} S \hspace{.15cm} $};
    \node at (0,1.6) (top2) [ellipse,draw]{$ \hspace{.15 cm} S \hspace{.15cm} $};
        \node at (top2.0) (invisible) [right] {};
    \node at (bottom2.0) (invisible2) [right] {};
    \filldraw (top2.180) circle (.5mm);
    \filldraw (bottom2.180) circle (.5mm);    

    \draw (JW.35) arc (180:0:.2cm) -- ++(-90:.7cm) arc (0:-180:.2cm) --(JW.-35);
    \node at (JW.75) [above] {...};
        \node at (JW.-75) [below] {...};
    \draw (JW.135) arc (180:0:1cm) -- ++(-90:.7cm) arc (0:-180:1cm) --(JW.-135);
        \draw (invisible.90)--++(90:1cm);
    \draw (invisible.90) .. controls ++(-90:6mm) and ++(90:10mm) .. (JW.145);
    \draw (invisible2.-90)--++(-90:1cm);
    \draw (invisible2.-90) .. controls ++(90:6mm) and ++(-90:10mm) .. (JW.-145);

    \draw (top2.35) -- ++(90:1cm);
    \draw (top2.60) -- ++(90:1cm);
    \node at (-.1,2.1) {...};
    \draw (top2.145) -- ++(90:1cm);

    \draw (bottom2.-35) -- ++(-90:1cm);
    \draw (bottom2.-60) -- ++(-90:1cm);
    \node at (-.1,-2.1) {...};
    \draw (bottom2.-145) -- ++(-90:1cm);
\end{tikzpicture}\\

%% file: diagrams/tikz/JWzerosimplify.tex
& =\left( \frac{1}{\qi{2n-1}} - \frac{\qi{2} }{\qi{2n-1} \qi{2}} \right)
\begin{tikzpicture}[baseline]
    \clip (-1,-1.3) rectangle (1,1.3);

    \node at (0,.6) (top) [ellipse,draw] {$ \hspace{.15 cm} S \hspace{.15cm} $};
    \node at (0,-.6) (bottom) [ellipse,draw]{$ \hspace{.15 cm} S \hspace{.15cm} $};
    \filldraw (top.180) circle (.5mm);
    \filldraw (bottom.180) circle (.5mm);    

    \node at (top.0) (invisible) [right] {};
    \draw (invisible.90)--++(90:2cm);
    \draw (invisible.90) -- ++ (-90:2cm);

    \draw (top.35) -- ++(90:1cm);
    \draw (top.60) -- ++(90:1cm);
    \node at (-.1,1.2) {...};
    \draw (top.145) -- ++(90:1cm);

    \draw (bottom.-35) -- ++(-90:1cm);
    \draw (bottom.-60) -- ++(-90:1cm);
    \node at (-.1,-1.2) {...};
    \draw (bottom.-145) -- ++(-90:1cm);
\end{tikzpicture} \\
& =0.

%% file: text/braiding.tex
\subsection{A partial braiding} \label{sec:braiding}

Recall the definition of a crossing given in \S \ref{ssTL}.  This still defines an element of $\pa$ and, away from $S$ boxes, diagrams related by a framed three-dimensional isotopy are equal in the planar algebra.  However, one needs to be careful manipulating these crossings and $S$ boxes at the same time.

\begin{thm}\label{thm:passacrossS}
You can isotope a strand above an $S$ box, but isotoping a strand below
an $S$ box introduces a factor of $-1$.
\begin{enumerate}
\item\label{pullover}
$\input{diagrams/tikz/pullstringoverS.tex}$

\item\label{pullunder}
$\input{diagrams/tikz/pullstringunderS.tex}$

\end{enumerate}
\end{thm}

\begin{proof}

\begin{enumerate}
\item
This is a straightforward consequence of part \ref{enoughcappings} of
Theorem \ref{easyconsequences}. Capping at any position from $2$ to
$4n-3$ in either of the above pictures gives zero; the first capping of
the lefthand picture is

\begin{align*}
\begin{tikzpicture}[baseline]
    \clip (-1.1,-1.4) rectangle (1.1,1.4);
    \node (S) at (0,0) [circle, draw] {$S$};
    \filldraw (S.180) circle (.5mm);
    \draw[rounded corners=2mm] (S.135) -- ++(90:.5cm) -- ++(180:.5cm) -- ++ (-90:1.7cm) -- ++(0:1.6cm) -- ++(90:3cm);
    \draw (S.90) node [above] {...};
    \draw (S.45) -- ++(90:15mm);
\end{tikzpicture} & =i S, \\
\intertext{and the first capping of the righthand picture is}
\input{diagrams/tikz/cappullstringoverS1.tex}
& = iq^{3/2} \input{diagrams/tikz/cappullstringoverS2.tex} \\
& = iq^{3/2} (iq^{1/2})^{4n-5} S \\
& = i S.
\end{align*}
Thus these two pictures are equal.

\item This is essentially identical to the previous argument, except that the factor picked up by resolving the crossings of the second picture is $$-iq^{-3/2} (-iq^{-1/2})^{4n-5} =- i;$$ hence the minus sign in the relation.
\end{enumerate}
\end{proof}

\begin{rem}
Upon reading this paper, one might hope that all subfactor planar algebras are braided, or partially braided.  Unfortunately this is far from being the case.  For the representation theory of the annular Temperley-Lieb category for $\qi{2} >2$, set out in \cite{MR1659204} and in the language of planar algebras in \cite{MR1929335}, implies that one cannot pull strands across lowest weight generators, even up to a multiple.  To see this, resolve all crossings in either of the equations in Theorem \ref{thm:passacrossS}; such an identity would give a linear dependence between ``annular consequences" of the generator.  For the other $\qi{2} < 2$ examples, namely $E_6$ and $E_8$,  \cite{math/0903.0144} shows that Equation \eqref{pullover} holds, but not Equation \eqref{pullunder}, even up to a coefficient.  The $\qi{2} = 2$ cases remain interesting.

\end{rem}

\begin{cor}
\label{cor:atmost1S}
Any diagram in $\pa$ is equal to a sum of diagrams involving at most one
$S$.
\end{cor}

\begin{proof}
When a diagram has more than one $S$, use the above relations to move one
of the $S$'s next to another one, then apply relation (\ref{twoS}) of
Definition \ref{def:pa} to replace the two $S$'s with a Jones-Wenzl
idempotent.  Resolve all the crossings and proceed inductively.
\end{proof}

\begin{cor}
\label{cor:evaluation} Every closed diagram is a multiple of the empty
diagram.
\end{cor}
\begin{proof}
By the previous corollary, a closed diagram can be written in terms of
closed diagrams with at most one $S$. If a closed diagram has exactly one
$S$, it must be zero, because the $S$ must have a cap attached to it
somewhere. If a closed diagram has no $S$'s, it can be rewritten as a
multiple of the empty diagram using Relation \ref{delta}, which allows us
to remove closed loops.
\end{proof}

\begin{cor}
\label{cor:spherical} The planar algebra $\pa$ is spherical.
\end{cor}
\begin{proof}
A braiding always suffices to show that a planar algebra is spherical; even
though there are signs when a strand passes underneath an $S$, we can
check that $\pa$ is spherical simply by passing strands above everything
else in the diagram.
\begin{equation*}
\mathfig{0.1}{spherical/1} = \mathfig{0.1}{spherical/2} =
\mathfig{0.1}{spherical/3}
\end{equation*}
\end{proof}

%
%

%% file: diagrams/tikz/pullstringoverS.tex
\begin{tikzpicture}[baseline]
    \clip (-1.1,-1.4) rectangle (1.1,1.4);

    \node (S) at (0,0) [circle, draw] {$S$};
     \filldraw (S.180) circle (.5mm);

    \draw[rounded corners=2mm] (-1,2) -- (-1,-1) -- (1,-1) -- (1,2);

    \draw (S.135) -- ++(90:15mm);
    \draw (S.90) node [above] {...};
    \draw (S.45) -- ++(90:15mm);
\end{tikzpicture}
=
\begin{tikzpicture}[baseline]
    \clip (-1.1,-1.4) rectangle (1.1,1.4);

    \node (S) at (0,0) [circle, draw] {$S$};
     \filldraw (S.180) circle (.5mm);

    \node (x1) at (S.135 |- -1,1){};
    \node (x2) at  (S.45 |- -1,1){} ;

    \draw[rounded corners=2mm] (-1,1.5) --  (-1,1) -- (1,1) -- (1,1.5);

    \draw (S.135) -- (x1.-90);
    \draw (x1.90) -- ++(90:5mm);
    \draw (S.90) node [above] {...};
    \draw (x2.90) -- ++(90:5mm);
    \draw (S.45) -- (x2.-90);
\end{tikzpicture}

%% file: diagrams/tikz/pullstringunderS.tex
\begin{tikzpicture}[baseline]
    \clip (-1.1,-1.4) rectangle (1.1,1.4);

    \node (S) at (0,0) [circle, draw] {$ S$};
     \filldraw (S.180) circle (.5mm);

    \draw[rounded corners=2mm] (-1,2) -- (-1,-1) -- (1,-1) -- (1,2);

    \draw (S.135) -- ++(90:15mm);
    \draw (S.90) node [above] {...};
    \draw (S.45) -- ++(90:15mm);
\end{tikzpicture}
=-
\begin{tikzpicture}[baseline]
    \clip (-1.1,-1.4) rectangle (1.1,1.4);

    \node (S) at (0,0) [circle, draw] {$S$};
     \filldraw (S.180) circle (.5mm);

    \node (x1) at (S.135 |- -1,1){};
    \node (x2) at  (S.45 |- -1,1){} ;

    \draw[rounded corners=2mm](-1,1.5) -- (-1,1) -- (x1.180);
    \draw (x1.0) -- (x2.180);
    \draw[rounded corners=2mm] (x2.0) -- (1,1)--(1,1.5);

    \draw (S.135) -- (x1.90);
    \draw (x1.90) -- ++(90:5mm);
    \draw (S.90) node [above] {...};
    \draw (x2.90) -- ++(90:5mm);
    \draw (S.45) -- (x2.90);
\end{tikzpicture}

%% file: diagrams/tikz/cappullstringoverS1.tex
\begin{tikzpicture}[baseline]
    \node (S) at (0,0) [circle, draw] {$S$};
     \filldraw (S.180) circle (.5mm);
    \node (x1) at (S.135 |- -1,1){};
    \node (x2) at  (S.45 |- -1,1){} ;
    \draw (S.135) -- (x1.-90);
    \draw[rounded corners=2mm] (x1.90) -- ++(90:5mm) -- (-1,1.65)--(-1,1) -- (1,1) -- (1,1.5);
    \draw (S.90) node [above] {...};
    \draw (x2.90) -- ++(90:5mm);
    \draw (S.45) -- (x2.-90);
\end{tikzpicture}

%% file: diagrams/tikz/cappullstringoverS2.tex
\begin{tikzpicture}[baseline]
    \node (S) at (0,0) [circle, draw] {$S$};
     \filldraw (S.180) circle (.5mm);
    \node (x1) at (S.135 |- -1,1){};
    \node (x2) at  (S.45 |- -1,1){} ;
    \draw[rounded corners=1mm] (S.135) -- (x1.-90) -- (x1.0);
    \draw[rounded corners=2mm] (x1.0) -- (x2.180) -- (1,1) -- (1,1.5);
    \draw (S.90) node [above] {...};
    \draw[rounded corners=2mm] (S.45) -- (x2.-90);
    \draw (x2.90) -- ++(90:5mm);
\end{tikzpicture}

%% file: text/consistency.tex
\subsection{The planar algebra $\pa$ is non-zero}
\label{sec:consistency}
In this section, we prove the following reassuring result.
\begin{thm}%
\label{thm:consistency}%
In the planar algebra $\pa$ described in
Definition \ref{def:pa}, the empty diagram is not equal to zero.
\end{thm}

The proof is fairly straightforward. We describe an algorithm for
evaluating any closed diagram in $\pa$, producing a number. Trivially,
the algorithm evaluates the empty diagram to $1$. We show that modifying
a closed diagram by one of the generating relations does not change the
result of the evaluation algorithm.

The algorithm we'll use actually allows quite a few choices along the
way, and the hard work will all be in showing that the answer we get
doesn't depend on these choices.\footnote{See \cite{MR1403861} and
\cite{MR2308953} for a related idea, sometimes called `confluence'.}
After that, checking that using a relation does not change the result
will be easy.

\begin{figure}[!ht]
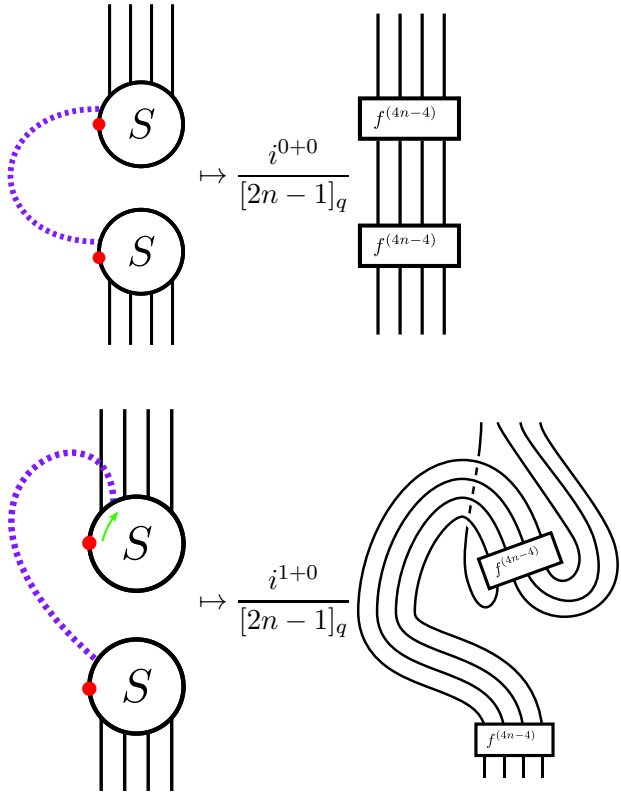

\begin{align*}
\mathfig{0.17}{consistency/arc} & \mapsto \frac{i^{0+0}}{\qi{2n-1}} \mathfig{0.1}{consistency/arc-connection} \\ \\
\mathfig{0.17}{consistency/arc2} & \mapsto \frac{i^{1+0}}{\qi{2n-1}} \mathfig{0.25}{consistency/arc-connection2}
\end{align*}
\caption{The main step of the algorithm chooses an arc connecting a pair of $S$ boxes, and replaces them with Jones-Wenzl idempotents connected by $4n-4$ parallel strands that cross above the original strands of the diagram,
and inserts a factor of $\frac{1}{\qi{2n-1}}$ and a power of $i$ determined by the attachment points of the arc.}
\label{fig:algorithm}
\end{figure}

\begin{defn}[Evaluation algorithm]
\label{defn:algorithm}
This is a function from closed diagrams (no relations) to the complex
numbers.
\begin{enumerate}
\item \label{step1}%
If there are at least two $S$ boxes in the diagram, choose a pair of $S$
boxes, and an imaginary arc connecting them. (The arc should be
transverse to everything in the diagram.)

Multiply by $\frac 1 {\qi{2n-1}}$, and replace the chosen pair of $S$ boxes with a pair of $f^{(4n-4)}$ boxes, connected by $4n-4$ parallel
strands following the arc which cross above any strands of the diagram
that the arc crosses (as illustrated in Figure \ref{fig:algorithm}).
Further, for each $S$ box do the following.  Starting at the marked point walk clockwise around the box counting the number of strands you pass before you reach the point where the arc attaches. This gives two numbers; multiply the new picture by $i$ raised to the sum of these two numbers.  Restart the algorithm on the result.
\item \label{step2}%
If there is exactly one $S$ box in the diagram, evaluate as $0$.
\item \label{step3}%
If there are no $S$ boxes in the diagram, evaluate the diagram in
Temperley-Lieb at $q + q^{-1} = 2 \cos(\frac{\pi}{4n-2})$.
\end{enumerate}
\end{defn}

\begin{thm}
\label{thm:algorithm-well-defined}
The algorithm is well-defined, and doesn't depend on any of the choices
made.
\end{thm}
\begin{proof}
We'll prove this in five stages.
\begin{enumerate}[i]
\item \label{alg:1}%
If two applications of the algorithm use the same pairing of $S$
boxes, and the same arcs, but replace the pairs in different orders, we
get the same answer.
\item \label{alg:2}%
If we apply the algorithm to a diagram with exactly two $S$ boxes,
then we can isotope the arc connecting them without affecting the answer.
\item \label{alg:3}%
Isotoping any arc does not change the answer.
\item \label{alg:4}%
Changing the point at which an arc attaches to an $S$ box does not
change the answer.
\item \label{alg:5}%
Two applications of the algorithm which use different pairings of the $S$
boxes give the same answers.
\end{enumerate}

\begin{proof}[Stage \ref{alg:1}]
Switching the order of two pairs of $S$ boxes produces Temperley-Lieb
diagrams that differ only where the corresponding two arcs intersected;
there we see one set of $4n-4$ parallel strands passing either over or
under the other set of $4n-4$ parallel strands. However, by Lemma
\ref{lem:f-crossings}, these are the same in Temperley-Lieb at $q + q^{-1} =
2 \cos(\frac{\pi}{4n-2})$. Thus after Step \ref{step3} of the evaluation
algorithm we get the same result.
\noqed\end{proof}
\begin{proof}[Stage \ref{alg:2}]
This follows easily from the fact that Temperley-Lieb is braided, and the final statement in
Lemma \ref{lem:f-twist}.
\noqed\end{proof}
\begin{proof}[Stage \ref{alg:3}]
In order to isotope an arbitrary arc, we make use of Stage \ref{alg:1} to
arrange that this arc corresponds to the final pair of $S$ boxes chosen.
Stage \ref{alg:2} then allows us to move the arc.
\noqed\end{proof}
\begin{proof}[Stage \ref{alg:4}]
Changing the point of attachment of an arc by one step clockwise results in a
Temperley-Lieb diagram at Step \ref{step3} which differs just by a factor
of $i$, according to the first part of Lemma \ref{lem:f-twist}. See the second part of Figure \ref{fig:algorithm}, which illustrates exactly this situation. This
exactly cancels with the factor of $i$ put in by hand by the algorithm.  Furthermore, moving the point of attachment across the marked point does not change the diagram, but does multiply it by a factor of
$i^{4n-4}=1$.
\noqed\end{proof}
\begin{proof}[Stage \ref{alg:5}]
We induct on the number of $S$ boxes in the diagram. If there are fewer
than $3$ $S$ boxes, there is no choice in the pairing. If there are
exactly $3$ $S$ boxes, the evaluation is automatically $0$.

Otherwise, consider two possible first choices of a pair of $S$ boxes.
Suppose one choice involves boxes which we'll call $A$ and $B$, while the
other involves boxes $C$ and $D$. There are two cases depending on
whether the sets $\{A,B\}$ and $\{C,D\}$ are disjoint, or have one common
element, say $D=A$. If the sets are disjoint, we (making use of the
inductive hypothesis), continue the algorithm which first removes $A$ and
$B$ by next removing $C$ and $D$, and continue the algorithm which first
removes $C$ and $D$ by next removing $A$ and $B$. The argument given in
Stage \ref{alg:1} shows that the final results are the same.
Alternatively, if the sets overlap, say with $A=D$, we choose some fourth $S$ box, say
$E$. After removing $A$ and $B$, we remove $C$ and $E$, while after
removing $A$ and $C$ we remove $B$ and $E$, and in each case we then
finish the algorithm making the same choices in either application. The
resulting Temperley-Lieb diagrams which we finally evaluate in Step
\ref{step3} differ exactly by the two sides of the identity in Lemma
\ref{lem:2f}. (More accurately, in the case that the arcs connecting
these pairs of $S$ boxes cross strands in the original diagram, the
resulting Temperley-Lieb diagrams differ by the two sides of that
equation sandwiched between some fixed Temperley-Lieb diagram; see
the remark following Lemma \ref{lem:2f}.)
\noqed\end{proof}
\end{proof}

\begin{proof}[Proof of Theorem \ref{thm:consistency}]
We just need to check that modifying a closed diagram by one of the
relations from Definition \ref{def:pa} does not change the answer.

\paragraph{Relation \eqref{delta}}
Make some set of choices for running the
algorithm, choosing arcs that avoid the disc in which the relation is
being applied. The set of choices is trivially valid both before and
after applying the relation. Once we reach Step \ref{step3} of the
algorithm, the Temperley-Lieb diagrams differ only by the relation, which
we know holds in Temperley-Lieb!

\paragraph{Relation \eqref{rotateS}}%
Run the algorithm, choosing the $S$ we want to rotate as one of the first
pair of $S$ boxes, using the same arc both before and after rotating the
$S$. The algorithm gives answers differing just by a factor of $i$,
agreeing with the relation. See Figure \ref{fig:rotateS}.

\paragraph{Relation \eqref{capS}}%
If there's exactly one $S$ box,
the algorithm gives zero anyway. If there's at least two $S$ boxes, choose the $S$ with a cap on it as a member of one of the pairs. Once
we reach Step \ref{step3}, the $S$ with a cap on it will have been
replaced with an $f^{(4n-4)}$ with a cap on one end, which gives $0$ in
Temperley-Lieb.

\paragraph{Relation \eqref{twoS}}%
When running the algorithm, on the diagram with more $S$ boxes, ensure that the pair of $S$ boxes affected
by the relation are chosen as a pair in Step \ref{step1}, with an arc
compatible with the desired application of the relation.

\begin{figure}[!ht]
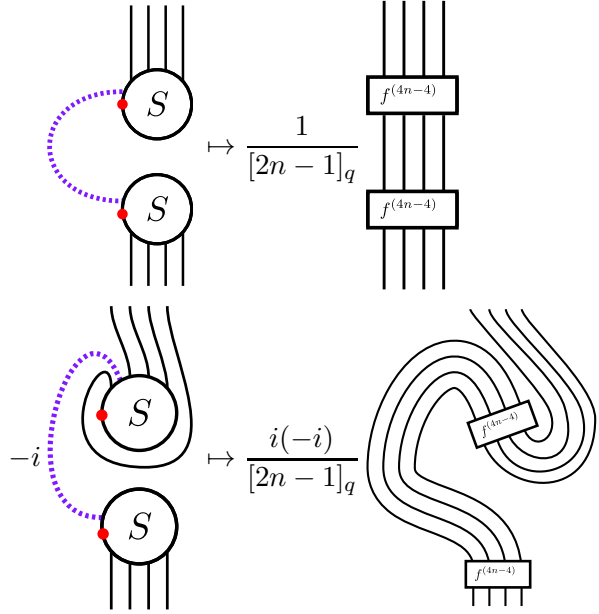

\begin{align*}
\mathfig{0.14}{consistency/arc} & \mapsto  \frac{1}{\qi{2n-1}} \mathfig{0.09}{consistency/arc-connection} \\
-i \mathfig{0.14}{consistency/rotation} & \mapsto  \frac{i(-i)}{\qi{2n-1}} \mathfig{0.22}{consistency/rotation-connection}
\end{align*}
\caption{The left hand sides are related by the rotation relation \ref{rotateS}. The algorithm gives the same answer in both cases.}
\label{fig:rotateS}
\end{figure}

\end{proof}

%% file: text/category.tex
\subsection{The tensor category of projections of a planar algebra}
\label{sec:category}%

In this section we describe a tensor category associated to a planar
algebra, whose objects are the `projections'. This is essentially
parallel to the construction of the tensor category of bimodules over a
subfactor \cite{MR1424954}. The tensor category described here is in fact
isomorphic to that one, although we won't need to make use of this fact.
We describe the category independently here, to emphasize that it can be
constructed directly from the planar algebra, without reference to the
associated subfactor.

\begin{defn}
Given a planar algebra $P$ we construct a tensor category $\cC_P$ as follows.
\begin{itemize}
\item An object of $\cC_P$ is a projection in one of the $2n$-box algebras $P_{2n}$.
\item Given two projections $\pi_1 \in P_{2n}$ and $\pi_2 \in P_{2m}$ we define $\Hom{}{\pi_1}{\pi_2}$ to be the space $\pi_2 P_{n \rightarrow m} \pi_1$  ($P_{n \rightarrow m}$ is a convenient way of denoting $P_{n+m}$, drawn with $n$ strands going down and $m$ going up.)
\item The tensor product $\pi_1 \otimes \pi_2$ is the disjoint union of the two projections in $P_{2n+2m}$.
\item The trivial object $\id$ is the empty picture (which is a projection in $P_0$).
\item The dual $\overline{\pi}$ of a projection $\pi$ is given by rotating it 180 degrees.
\end{itemize}
\end{defn}

This category comes with a special object $X \in P_2$ which is the single strand.  Note that $X = \overline{X}$. 
We would like to be able to take direct sums in this category.  If $\pi_1$ and $\pi_2$ are orthogonal projections in the same box space $P_n$ (i.e. if $\pi_1 \pi_2 = 0 = \pi_2 \pi_1$), then their direct sum is just $\pi_1 + \pi_2$.  However, if they are not orthogonal the situation is a bit more difficult.  One solution to this problem is to replace the projections with isomorphic projections which are orthogonal.  However, this construction only makes sense on equivalence classes, so we use another construction.

\begin{defn}
Given a category $\cC$ define its matrix category $\Mat{\cC}$ as follows.
\begin{enumerate}
\item The objects of $\Mat{\cC}$ are formal direct sums of objects of $\cC$
\item A morphism of $\Mat{\cC}$ from $A_1 \oplus \ldots \oplus A_n \rightarrow B_1 \oplus \ldots \oplus B_m$ is an $m$-by-$n$ matrix whose $(i,j)$th entry is in $\Hom{\cC}{A_j}{B_i}$.
\end{enumerate}
\end{defn}


If $\cC$ is a tensor category then $\Mat{\cC}$ has an obvious tensor
product (on objects, formally distribute, and on morphisms, use the
usual tensor product of matrices and the tensor product for $\cC$ on
matrix entries).  If $\cC$ is spherical then so is $\Mat{\cC}$ (where the
dual on objects is just the dual of each summand and on morphisms the dual transposes
the matrix and dualizes each matrix entry).

\begin{lem}
If $\pi_1$ and $\pi_2$ are orthogonal projections in $P_{2n}$, then $\pi_1 \oplus \pi_2 \cong \pi_1 + \pi_2$ in $\Mat{\cC_P}$.
\end{lem}
\begin{proof}
Define $f: \pi_1 \oplus \pi_2 \rightarrow \pi_1 + \pi_2$ and $g: \pi_1 + \pi_2 \rightarrow \pi_1 \oplus \pi_2$ by $f = \begin{pmatrix} \pi_1 & \pi_2 \end{pmatrix}$ and $g=\begin{pmatrix} \pi_1 \\ \pi_2 \end{pmatrix}$.  Then $f g = \id_{\pi_1 + \pi_2}$ and $g f= \id_{\pi_1 \oplus \pi_2}$.
\end{proof}

\begin{defn}
An object $\pi \in \Mat{\cC_P}$ is called minimal if $\Hom{}{\pi}{\pi}$ is $1$-dimensional.
\end{defn}
\begin{defn}
A planar algebra is called semisimple if every projection is a direct sum of minimal projections and for any pair of non-isomorphic minimal projections $\pi_1$ and $\pi_2$, we have that $\Hom{}{\pi_1}{ \pi_2} = 0$.
\end{defn}

\begin{defn}
The principal (multi-)graph of a semisimple planar algebra has as vertices the isomorphism classes of minimal projections, and there are $$\dim \Hom{}{\pi_1 \otimes X}{\pi_2} (= \dim \Hom{}{\pi_1}{\pi_2 \otimes X})$$ edges between the vertices $\pi_1 \in P_n$ and $\pi_2 \in P_m$.
\end{defn}

Our definitions here are particularly simple because we work in the
context of unshaded planar algebras. A slight variation works for a
shaded planar algebra as well. \footnote{
For a shaded planar algebra, we construct a $2$-category with two
objects `shaded' and `unshaded'. The $1$-morphisms from $A$ to $B$
are now the projections in the planar algebra, where $A$ is the shading at the marked point, and $B$ is the shading opposite the marked point. The $2$-morphisms are defined the same way as the morphisms were above. This $2$-category is equivalent to the $2$-category of $N,M$ bimodules for the type $II_1$ subfactor $N \subset M$ associated to the shaded planar algebra, as described in \cite{MR1424954}. }

\subsubsection{The minimal projections of $\pa$}

Now we can use the definitions of the previous section to explain why we call $\pa$ the $D_{2n}$ planar algebra.

\begin{thm}
The planar algebra $D_{2n}$ is semi-simple, with minimal projections $f^{(k)}$ for $k=0,
\ldots, 2n-3$ along with $P$ and $Q$ defined by
\begin{align*}
P & = \frac{1}{2}\left(f^{(2n-2)} + S\right) \\
\intertext{and}
Q & = \frac{1}{2}\left(f^{(2n-2)} - S\right).
\end{align*}
The principal graph is the Dynkin diagram $D_{2n}$.
\begin{equation*}
\mathfig{0.5}{graphs/d2n}
\end{equation*}
\end{thm}
\begin{proof} Observe that $\JW{2n-2} \cdot S = S$ (as the identity has weight $1$ in $\JW{2n-2}$ and all non-identity Temperley-Lieb pictures have product $0$ with $S$) and $S^2=\JW{2n-2}$. We see that $P$ and $Q$ are projections.
Let $\mathcal{M} = \{f^{(1)}, \ldots f^{(2n-3)}, P, Q\}$.  By Lemmas \ref{LemmasA1} and \ref{LemmasA2} (below) every projection in $\mathcal{M}$ is minimal. Lemma \ref{lem:no-homs} says there are no nonzero morphisms between different elements of $\mathcal{M}$. By Lemmas \ref{LemmasB1}, \ref{LemmasB2}, and \ref{LemmasB3}, we see that for each $Y \in \mathcal{M}$, the projection $Y \tensor f^{(1)}$ is isomorphic to a direct sum of projections in $\mathcal{M}$.  Thus, because every projection is a summand of $\id \in P_n$ for some $n$, every minimal projection is in $\mathcal{M}$.  Finally from Lemmas \ref{LemmasB1}, \ref{LemmasB2}, and \ref{LemmasB3} we read off that the principal graph for our planar algebra is the Dynkin diagram $D_{2n}$.
\end{proof}

\begin{rem}
Since $S^* = S$, all the projections are self-adjoint. The projections
$f^{(k)}$ are all self-dual. The projections $P$ and $Q$ are self-dual
when $n$ is odd, and when $n$ is even, $\overline{P} = Q$ and
$\overline{Q} = P$. These facts follow immediately from the definitions,
and the rotation relation (\ref{rotateS}).
\end{rem}

\begin{lem} \label{LemmasA1}
The Jones-Wenzl idempotents $f^{(k)}$ for $k=0, \ldots, f^{(2n-3)}$ are minimal.
\end{lem}

\begin{rem}
The minimality of the empty diagram, $f^{(0)}$, is exactly the fact that any closed diagram evaluates to a multiple of the empty diagram; that is, $\dim{\pa_0}=1$.
\end{rem}

\begin{proof}
The space $\Hom{}{f^{(i)}} {f^{(i)}}$ consists of all diagrams obtained by filling in the empty ellipse in the following diagram.

\begin{equation*}
\mathfig{0.14}{projections/fi}
\end{equation*}

We want to show that any such diagram which is non-zero is equal to a multiple of the diagram gotten by inserting the identity into the empty ellipse. By Corollary \ref{cor:atmost1S}, we need only consider diagrams with $0$ or $1$ $S$ boxes. First consider inserting any Temperley-Lieb diagram.  Since any cap applied to a Jones-Wenzl is zero, the ellipse must contain no cups or caps, hence it is a multiple of the identity.  Now consider any diagram with exactly one $S$.  Since $S$ has $4n-4$ strands, and $2i \leq 4n-6$, any such diagram must cap off the $S$, hence it vanishes.
\end{proof}

\begin{lem} \label{LemmasA2}
The projections $P = \frac{1}{2}\left(f^{(2n-2)} + S\right)$ and $Q = \frac{1}{2}\left(f^{(2n-2)} - S\right)$ are minimal.
\end{lem}
\begin{proof}
The two proofs are identical, so we do the $P$ case.  The space $\Hom{}{P} {P}$ consists of all ways of filling in the following diagram.

\begin{equation*}
\mathfig{0.13}{projections/PP}
\end{equation*}

We want to show that any such diagram which is non-zero is equal to a multiple of the diagram with the identity inserted.  Again we use Corollary \ref{cor:atmost1S}. First consider any Temperley-Lieb diagram drawn there. Since any cap applied to $P$ is zero, the diagram must have no cups or caps, hence it is a multiple of the identity.  Now consider any diagram with exactly one $S$.  Since $S$ has $4n-4$ strands, any such diagram which does not cap off $S$ must be (up to rotation) the following diagram.

\begin{equation*}
\mathfig{0.13}{projections/PPS}
\end{equation*}

Since $PS=P$, this diagram is a multiple of the diagram with the identity inserted.

\end{proof}

\begin{lem}
\label{lem:no-homs}
If $A$ and $B$ are two distinct projections from the set $$\{f^{(0)}, f^{(1)}, \ldots, f^{(2n-3)}, P, Q\}$$ then
$\Hom{}{A}{B} = 0$.
\end{lem}
\begin{proof}
Suppose $A$ and $B$ are distinct Jones-Wenzl projections.  Any morphism between them with exactly one $S$ must cap off the $S$, and so is $0$.  Any morphism between them in Temperley-Lieb must cap off either $A$ or $B$ and so is zero.  If $A$ is a Jones-Wenzl projection, while $B$ is $P$ or $Q$, exactly the same argument holds.

If $A=P$ and $B=Q$, we see that the morphism space is spanned by Temperley-Lieb diagrams and the diagram with a single $S$ box. Changing basis, the morphism space is spanned by non-identity Temperley-Lieb diagrams, along with $P$ and $Q$. Non-identity Temperley-Lieb diagrams are all zero as morphisms, because they result in attaching a cap to both $P$ and $Q$. The elements $P$ and $Q$ are themselves zero as morphisms from $P$ to $Q$, because $PQ=QP=0$.
\end{proof}

\begin{lem} \label{LemmasB1}
The projection $f^{(k)} \tensor f^{(1)}$ is isomorphic to $f^{(k-1)} \directSum f^{(k+1)}$ for $k = 1,\ldots, 2n-4$.
\end{lem}
\begin{proof}
This is a well known result about Temperley-Lieb.   The explicit isomorphisms are

\begin{align*}
\begin{pmatrix}
\frac{\qi{m}}{\qi{m+1}} \cdot \mathfig{0.1}{projections/capk_fX_f} \\ \mathfig{0.1}{projections/id_fX_f}
\end{pmatrix}: &\mathfig{0.1}{projections/fk1} \to f^{(k-1)} \directSum f^{(k+1)} \\
\intertext{and}
\begin{pmatrix}
\mathfig{0.1}{projections/cupk_f_fX}  & \mathfig{0.1}{projections/id_f_fX}
\end{pmatrix}: &  f^{(k-1)} \directSum f^{(k+1)} \to \mathfig{0.1}{projections/fk1}
.\end{align*}
The fact that these are inverses to each other is exactly Wenzl's relation \eqref{WenzlsRelation}.
\end{proof}

\begin{lem} \label{LemmasB2}
The projection $f^{(2n-3)} \tensor f^{(1)}$ is isomorphic to $f^{(2n-4)} \directSum P \directSum Q$.
\end{lem}
\begin{proof}
The explicit isomorphisms are:

\begin{align*}
\begin{pmatrix}
\frac{\qi{m}}{\qi{m+1}} \cdot \mathfig{0.1}{projections/capn_fX_f} \\ \mathfig{0.1}{projections/id_fX_P} \\ \mathfig{0.1}{projections/id_fX_Q}
\end{pmatrix}: \mathfig{0.1}{projections/fn1} \to f^{(2n-4)} \directSum P \directSum Q \\
\intertext{and}
\begin{pmatrix}
\mathfig{0.1}{projections/cupn_f_fX} & \mathfig{0.1}{projections/id_P_fX} & \mathfig{0.1}{projections/id_Q_fX}
\end{pmatrix}: f^{(2n-4)} \directSum P \directSum Q \to \mathfig{0.1}{projections/fn1}.
\end{align*}
The fact that these are inverses to each other follows from Wenzl's relation and the fact that $P$ and $Q$ absorb Jones-Wenzl idempotents (ie, $f^{(2n-3)} \cdot P=P$ and $f^{(2n-3)} \cdot Q=Q$).
\end{proof}

\begin{lem} \label{LemmasB3}
$P \tensor f^{(1)} \iso f^{(2n-3)}$ and $Q \tensor f^{(1)} \iso f^{(2n-3)}$.
\end{lem}

\begin{proof}
We claim that the maps
$$\qi{2} \cdot
\begin{tikzpicture}[baseline]
    \node at (.2,.5) (mid) [rectangle, draw] {$ \hspace{.7cm} P \hspace{.7cm} $};
    \node at (0,-.5) (bottom) [rectangle,draw]{$ \JW{2n-3} $};
    \draw (mid.-40)--(bottom.35);
    \draw (mid.-20) arc (-180:0:.3cm) -- ++(90:.5cm);
    \draw (mid.-160)--(bottom.145);
     \node at (mid.-90) [below] {...};
\end{tikzpicture}
: P \tensor f^{(1)} \rightarrow f^{(2n-3)}$$
and
$$
\begin{tikzpicture}[baseline]
    \node at (0,.5) (top) [rectangle,draw] {$ \JW{2n-3} $};
    \node at (.2,-.5) (mid) [rectangle, draw] {$ \hspace{.7cm} P \hspace{.7cm} $};
    \draw (top.-35)--(mid.40);
    \draw (top.-145)--(mid.160);
    \draw (mid.20) arc (180:0:.3cm) -- ++(-90:.5cm);
    \node at (mid.90) [above] {...};
\end{tikzpicture}
:f^{(2n-3)} \rightarrow P \tensor f^{(1)}$$
are isomorphisms inverse to each other.  To check this, we need to verify
$$\qi{2} \cdot
\begin{tikzpicture}[baseline]
    \node at (0,1) (top) [rectangle,draw] {$ \JW{2n-3} $};
    \node at (.2,0) (mid) [rectangle, draw] {$ \hspace{.7cm} P \hspace{.7cm} $};
    \node at (0,-1) (bottom) [rectangle,draw]{$ \JW{2n-3} $};
    \draw (top.-35)--(mid.40);
    \draw (top.-145)--(mid.160);
    \draw (mid.-40)--(bottom.35);
    \draw (mid.20) arc (180:0:.3cm) -- ++(-90:.5cm);
    \draw (mid.-20) arc (-180:0:.3cm) -- ++(90:.5cm);
    \draw (mid.-160)--(bottom.145);
    \node at (mid.90) [above] {...};
     \node at (mid.-90) [below] {...};
\end{tikzpicture}
= \JW{2n-3} \text{ \qquad and \qquad }
\qi{2} \cdot
\begin{tikzpicture}[baseline]
    \node at (0,0) (top) [rectangle,draw] {$ \JW{2n-3} $};
    \node at (.2,-1) (mid) [rectangle, draw] {$ \hspace{.7cm} P \hspace{.7cm} $};

    \draw (top.-35)--(mid.40);
    \draw (top.-145)--(mid.160);
     \draw (mid.20) arc (180:0:.3cm) -- ++(-90:.5cm);

    \node at (mid.90) [above] {...};

    \node at (.2,1) (toptop) [rectangle, draw] {$ \hspace{.7cm} P \hspace{.7cm} $};
    \draw (top.35)--(toptop.-40);
    \draw (top.145)--(toptop.-160);
     \draw (toptop.-20) arc (-180:0:.3cm) -- ++(90:.5cm);

    \node at (toptop.-90) [below] {...};
\end{tikzpicture}
= P \tensor \JW{1}.$$

The first equality is straightforward:  capping $P=\frac{1}{2}(\JW{2n-2} + S)$ on the right side kills its $S$ component, and then the equality follows from the partial trace relation and the observation that  $\frac{\qi{2n-1}}{\qi{2n-2}}=\frac{2}{\qi{2}}$. The second will take a bit more work to establish, but it's not hard. We first observe that $P \cdot \JW{2n-3} = P$, then expand both $P$s as $\frac{1}{2} (\JW{2n-2} + S)$. Thus
\begin{align*}
\input{diagrams/tikz/PEnP.tex} \\
\input{diagrams/tikz/PEnPExpanded.tex},
\intertext{and applying Wenzl's relation to the first term and the two-S relation to the fourth term yields}
& = \input{diagrams/tikz/PEnPWenzl.tex}. \\
\intertext{Because $\qi{k}=\qi{4n-2-k}$ we can cancel the second and fifth terms; and using $\frac{\qi{2n-1}}{\qi{2n-2}}=\frac{2}{\qi{2}}$ again we get }
\input{diagrams/tikz/PEnPcancel.tex}.
\end{align*}

At this point, we use part \ref{enoughcappings} of Theorem \ref{easyconsequences} to show

$$\input{diagrams/tikz/JWEnS.tex}.$$

Start at the top left of the pictures, and take the first $4n-4$ cappings in the counterclockwise direction.  Most of these give zero immediately, and the three we are left to check are

\begin{align*}
\input{diagrams/tikz/JWEnStopcap.tex}, \\
\input{diagrams/tikz/JWEnSbottomcap.tex} \\
\intertext{and}
\input{diagrams/tikz/JWS.tex}.
\end{align*}

The first two of these follows from the partial trace relation (Equation \eqref{eq:PartialTraceRelation}) and $\JW{2n-3} \cdot S =S$, and the third follows from $\JW{2n-2} \cdot S =S$.  Therefore, we conclude
$$
\qi{2} \cdot
\begin{tikzpicture}[baseline]
    \node at (0,0) (top) [rectangle,draw] {$ \JW{2n-3} $};
    \node at (.2,-1) (mid) [rectangle, draw] {$ \hspace{.7cm} P \hspace{.7cm} $};

    \draw (top.-35)--(mid.40);
    \draw (top.-145)--(mid.160);
     \draw (mid.20) arc (180:0:.3cm) -- ++(-90:.5cm);

    \node at (mid.90) [above] {...};

    \node at (.2,1) (toptop) [rectangle, draw] {$ \hspace{.7cm} P \hspace{.7cm} $};
    \draw (top.35)--(toptop.-40);
    \draw (top.145)--(toptop.-160);
     \draw (toptop.-20) arc (-180:0:.3cm) -- ++(90:.5cm);

    \node at (toptop.-90) [below] {...};
\end{tikzpicture}
=
\frac{\qi{2}}{4} \cdot \frac{2}{\qi{2}}
\left(
\begin{tikzpicture}[baseline]
    \clip (-1,-1.3) rectangle (1,1.3);
    \node at (0,0) (box) [rectangle,draw] {$\JW{2n-2}$};

    \draw (box.-145)--++(-90:1cm);
    \node at (-.1,-.6) {...};
    \draw (box.-35) -- ++(-90:1cm);
    \draw (box.-60) -- ++(-90:1cm);

    \draw (box.35) -- ++(90:1cm);
    \draw (box.60) -- ++(90:1cm);
    \node at (-.1,.6) {...};
    \draw (box.145) -- ++(90:1cm);

    \draw (.9,1.5)--(.9,-1.5);
\end{tikzpicture}
+
\begin{tikzpicture}[baseline]
    \clip (-.8,-1.3) rectangle (1,1.3);
    \node at (0,0) (box) [ellipse,draw] {$ \hspace{.15 cm} S \hspace{.15cm} $};
    \filldraw (box.180) circle (.5mm);

    \draw (box.-145)--++(-90:1cm);
    \node at (-.1,-.6) {...};
    \draw (box.-35) -- ++(-90:1cm);

    \draw (box.35) -- ++(90:1cm);
    \node at (-.1,.6) {...};
    \draw (box.145) -- ++(90:1cm);

    \draw (.9,1.5)--(.9,-1.5);
\end{tikzpicture}
\right)
=
\begin{tikzpicture}[baseline]
    \clip (-1,-1.3) rectangle (1,1.3);
    \node at (0,0) (box) [rectangle,draw] {$\hspace{3mm} P \hspace{3mm}$};

    \draw (box.-145)--++(-90:1cm);
    \node at (-.1,-.6) {...};
    \draw (box.-35) -- ++(-90:1cm);
    \draw (box.-60) -- ++(-90:1cm);

    \draw (box.35) -- ++(90:1cm);
    \draw (box.60) -- ++(90:1cm);
    \node at (-.1,.6) {...};
    \draw (box.145) -- ++(90:1cm);

    \draw (.9,1.5)--(.9,-1.5);
\end{tikzpicture},
$$
which is what we wanted to show.
\end{proof}

\newcommand{\Wfour}{\DirectSum_{l=0}^{\frac{n-3}{2}} f^{(4l)}}
\newcommand{\Wfourp}{\DirectSum_{l=0}^{\frac{n-4}{2}} f^{(4l+2)}}


\subsubsection{The tensor product decompositions} We do not prove the formulas that follow, and they are not essential to this paper.  Nevertheless, we include the full tensor product table of $D_{2n}$ for the sake of making this description of $D_{2n}$ as complete as possible. Partial
tensor product tables appears in \cite[\S 3.5]{MR1145672} and \cite[\S 7]{MR1936496}.

Using the methods of this paper, one could prove that these tensor product formulas hold by producing explicit bases for all the appropriate $\operatorname{Hom}$ spaces in the tensor category of projections.  However, this method would not show that these formulas are the only extension of the data encoded in the principal graph.

Much of this is proved in \cite{MR1145672}, except for the formula for $f^{(j)} \tensor f^{(k)}$ when $2n-2 \leq j+k \leq 4n-4$ in Equation \eqref{eq:wjwk-large}
and the formula for $P \tensor f^{(2k+1)}$ and $Q \tensor f^{(2k+1)}$ in Equation \eqref{eq:pqw-odd}. Nonetheless, the methods of \cite{MR1145672} readily extend to give the remaining formulas. With the same exceptions, along with Equations \eqref{eq:pw-even} and \eqref{eq:qw-even},
these are proved in \cite{MR1936496}, by quite different methods.
Further, \cite{MR1145672} proves there is no associative tensor product
extending the tensor product data encoded in the principal graphs $D_{2n+1 \geq 5}$ with an odd number of vertices.

\begin{thm}
The tensor product structure is commutative, and described by the following isomorphisms.

When $j+k < 2n-2$
\begin{equation*}
f^{(j)} \tensor f^{(k)}  \iso \DirectSum_{l=\frac{\abs{j-k}}{2}}^{\frac{j+k}{2}} f^{(2l)}
\end{equation*}
and when $2n-2 \leq j+k \leq 4n-4$
\begin{equation}
\label{eq:wjwk-large}
f^{(j)} \tensor f^{(k)}  \iso \begin{cases}
 \left(\displaystyle \DirectSum_{l=\frac{\abs{j-k}}{2}}^{2n-3-\frac{j+k}{2}} f^{(2l)}\right) \directSum \left(\displaystyle \DirectSum_{l=2n-2-\frac{j+k}{2}}^{n-2} 2 f^{(2l)} \right) \directSum P \directSum Q  & \text{if $j+k$ is even} \\
 \left(\displaystyle \DirectSum_{l=\frac{\abs{j-k}}{2}}^{2n-3-\frac{j+k}{2}} f^{(2l)}\right) \directSum \left(\displaystyle \DirectSum_{l=2n-2-\frac{j+k}{2}}^{n-\frac{3}{2}} 2 f^{(2l)} \right) & \text{if $j+k$ is odd.}
\end{cases}
\end{equation}
Moreover
\begin{equation}
\label{eq:pqw-odd}%
P \tensor f^{(2k+1)} \iso Q \tensor f^{(2k+1)}
                    \iso \DirectSum_{l=0}^k f^{(2n-2l-3)}
\end{equation}
\begin{align}
\label{eq:pw-even}
P \tensor f^{(2k)}     & \iso \begin{cases}
 P \directSum \DirectSum_{l=0}^{k-1} f^{(2n-2l-4)} & \text{if $k$ is even} \\
 Q \directSum \DirectSum_{l=0}^{k-1} f^{(2n-2l-4)} & \text{if $k$ is odd}
 \end{cases} \\
\label{eq:qw-even}
Q \tensor f^{(2k)}     & \iso \begin{cases}
 Q \directSum \DirectSum_{l=0}^{k-1} f^{(2n-2l-4)} & \text{if $k$ is even} \\
 P \directSum \DirectSum_{l=0}^{k-1} f^{(2n-2l-4)} & \text{if $k$ is odd}
 \end{cases} \\
\intertext{and when $n$ is even}
\notag P \tensor P & \iso Q \directSum \Wfourp   \\
\notag Q \tensor Q & \iso P \directSum \Wfourp  \\
\notag P \tensor Q & \iso Q \tensor P \iso \Wfour \\
\intertext{and when $n$ is odd}
\notag P \tensor P & \iso P \directSum \Wfour \\
\notag Q \tensor Q & \iso Q \directSum \Wfour \\
\notag P \tensor Q & \iso Q \tensor P \iso  \Wfourp.
\end{align}
\end{thm}

%% file: diagrams/tikz/PEnP.tex
\qi{2} \cdot
\begin{tikzpicture}[baseline]
    \clip (-.6,-1.3) rectangle (.8,1.3);
    \node at (0,.7) (top) [rectangle,draw] {$ \hspace{.3 cm} P \hspace{.3cm} $};
    \node at (0,-.7) (bottom) [rectangle,draw]{$ \hspace{.3cm} P \hspace{.3cm} $};

    \draw (top.-60)--(bottom.60);
    \draw (top.-145)--(bottom.145);
    \node at (-.1,0) {...};
    \draw (top.-35) arc (-180:0:.2cm) -- ++(90:1cm);
    \draw (bottom.35) arc (180:0:.2cm) -- ++(-90:1cm);

    \draw (top.35) -- ++(90:1cm);
    \draw (top.60) -- ++(90:1cm);
    \node at (-.1,1.2) {...};
    \draw (top.145) -- ++(90:1cm);

    \draw (bottom.-35) -- ++(-90:1cm);
    \draw (bottom.-60) -- ++(-90:1cm);
    \node at (-.1,-1.2) {...};
    \draw (bottom.-145) -- ++(-90:1cm);
\end{tikzpicture}
& = \qi{2} \cdot
\begin{tikzpicture}[baseline]
    \clip (-1.5,-1.3) rectangle (2,1.3);
    \node at (0,.7) (top) [rectangle,draw] {$\frac{1}{2}(\JW{2n-2}+S)$};
    \node at (0,-.7) (bottom) [rectangle,draw]{$\frac{1}{2}(\JW{2n-2}+S)$};

    \draw (top.-30)--(bottom.30);
    \draw (top.-163)--(bottom.163);
    \node at (-.1,0) {...};
    \draw (top.-17) arc (-180:0:.2cm) -- ++(90:1cm);
    \draw (bottom.17) arc (180:0:.2cm) -- ++(-90:1cm);

    \draw (top.17) -- ++(90:1cm);
    \draw (top.30) -- ++(90:1cm);
    \node at (-.1,1.2) {...};
    \draw (top.163) -- ++(90:1cm);

    \draw (bottom.-17) -- ++(-90:1cm);
    \draw (bottom.-30) -- ++(-90:1cm);
    \node at (-.1,-1.2) {...};
    \draw (bottom.-163) -- ++(-90:1cm);
\end{tikzpicture}

%% file: diagrams/tikz/PEnPExpanded.tex
& = \frac{ \qi{2}}{4} \cdot
 \left(
\begin{tikzpicture}[baseline]
    \clip (-.8,-1.3) rectangle (1,1.3);
    \node at (0,.7) (top) [rectangle,draw] {$  \JW{2n-2} $};
    \node at (0,-.7) (bottom) [rectangle,draw] {$  \JW{2n-2} $};

    \draw (top.-60)--(bottom.60);
    \draw (top.-145)--(bottom.145);
    \node at (-.1,0) {...};
    \draw (top.-35) arc (-180:0:.2cm) -- ++(90:1cm);
    \draw (bottom.35) arc (180:0:.2cm) -- ++(-90:1cm);

    \draw (top.35) -- ++(90:1cm);
    \draw (top.60) -- ++(90:1cm);
    \node at (-.1,1.2) {...};
    \draw (top.145) -- ++(90:1cm);

    \draw (bottom.-35) -- ++(-90:1cm);
    \draw (bottom.-60) -- ++(-90:1cm);
    \node at (-.1,-1.2) {...};
    \draw (bottom.-145) -- ++(-90:1cm);
\end{tikzpicture}
+
\begin{tikzpicture}[baseline]
    \clip (-.8,-1.3) rectangle (1,1.3);

    \node at (0,.7) (top) [rectangle,draw] {$  \JW{2n-2} $};
    \node at (0,-.7) (bottom) [ellipse,draw]{$ \hspace{.15 cm} S \hspace{.15cm} $};
    \filldraw (bottom.180) circle (.5mm);    

    \draw (top.-60)--(bottom.60);
    \draw (top.-145)--(bottom.145);
    \node at (-.1,0) {...};
    \draw (top.-35) arc (-180:0:.2cm) -- ++(90:1cm);
    \draw (bottom.35) arc (180:0:.2cm) -- ++(-90:1cm);

    \draw (top.35) -- ++(90:1cm);
    \draw (top.60) -- ++(90:1cm);
    \node at (-.1,1.2) {...};
    \draw (top.145) -- ++(90:1cm);

    \draw (bottom.-35) -- ++(-90:1cm);
    \draw (bottom.-60) -- ++(-90:1cm);
    \node at (-.1,-1.2) {...};
    \draw (bottom.-145) -- ++(-90:1cm);
\end{tikzpicture}
+
\begin{tikzpicture}[baseline]
    \clip (-.8,-1.3) rectangle (1,1.3);

    \node at (0,.7) (top) [ellipse,draw] {$ \hspace{.15 cm} S \hspace{.15cm} $};
    \node at (0,-.7) (bottom) [rectangle,draw] {$  \JW{2n-2} $};
    \filldraw (top.180) circle (.5mm);

    \draw (top.-60)--(bottom.60);
    \draw (top.-145)--(bottom.145);
    \node at (-.1,0) {...};
    \draw (top.-35) arc (-180:0:.2cm) -- ++(90:1cm);
    \draw (bottom.35) arc (180:0:.2cm) -- ++(-90:1cm);

    \draw (top.35) -- ++(90:1cm);
    \draw (top.60) -- ++(90:1cm);
    \node at (-.1,1.2) {...};
    \draw (top.145) -- ++(90:1cm);

    \draw (bottom.-35) -- ++(-90:1cm);
    \draw (bottom.-60) -- ++(-90:1cm);
    \node at (-.1,-1.2) {...};
    \draw (bottom.-145) -- ++(-90:1cm);
\end{tikzpicture}
+
\begin{tikzpicture}[baseline]
    \clip (-.6,-1.3) rectangle (1,1.3);

    \node at (0,.6) (top) [ellipse,draw] {$ \hspace{.15 cm} S \hspace{.15cm} $};
    \node at (0,-.6) (bottom) [ellipse,draw]{$ \hspace{.15 cm} S \hspace{.15cm} $};
    \filldraw (top.180) circle (.5mm);
    \filldraw (bottom.180) circle (.5mm);    

    \draw (top.-60)--(bottom.60);
    \draw (top.-145)--(bottom.145);
    \node at (-.1,0) {...};
    \draw (top.-35) arc (-180:0:.2cm) -- ++(90:1cm);
    \draw (bottom.35) arc (180:0:.2cm) -- ++(-90:1cm);

    \draw (top.35) -- ++(90:1cm);
    \draw (top.60) -- ++(90:1cm);
    \node at (-.1,1.2) {...};
    \draw (top.145) -- ++(90:1cm);

    \draw (bottom.-35) -- ++(-90:1cm);
    \draw (bottom.-60) -- ++(-90:1cm);
    \node at (-.1,-1.2) {...};
    \draw (bottom.-145) -- ++(-90:1cm);
\end{tikzpicture}
\right)

%% file: diagrams/tikz/PEnPWenzl.tex
\frac{ \qi{2}}{4} \cdot
\left(
\left(
\frac{\qi{2n-1}}{\qi{2n-2}}
\begin{tikzpicture}[baseline]
    \clip (-1,-1.3) rectangle (1,1.3);
    \node at (0,0) (box) [rectangle,draw] {$\JW{2n-2}$};

    \draw (box.-145)--++(-90:1cm);
    \node at (-.1,-.6) {...};
    \draw (box.-35) -- ++(-90:1cm);
    \draw (box.-60) -- ++(-90:1cm);

    \draw (box.35) -- ++(90:1cm);
    \draw (box.60) -- ++(90:1cm);
    \node at (-.1,.6) {...};
    \draw (box.145) -- ++(90:1cm);

    \draw (.9,1.5)--(.9,-1.5);
\end{tikzpicture}
-
\frac{\qi{2n-1}}{\qi{2n-2}}
\begin{tikzpicture}[baseline]
    \clip (-1,-1.3) rectangle (1,1.3);
    \node at (0,0) (box) [rectangle,draw] {$\hspace{.1cm} \JW{2n-1} \hspace{.1cm} $};

    \draw (box.-50)-- ++(-90:1cm);
    \draw (box.-90) -- ++(-90:1cm);
    \draw (box.-150)--++(-90:1cm);
    \node at (-.25,-.6) {...};
    \draw (box.-30) -- ++(-90:1cm);

    \draw (box.50)-- ++(90:1cm);
    \draw (box.90) -- ++(90:1cm);
    \draw (box.150)--++(90:1cm);
    \node at (-.25,.6) {...};
    \draw (box.30) -- ++(90:1cm);
\end{tikzpicture}
\right) \right. \\
& \left. +
\begin{tikzpicture}[baseline]
    \clip (-.8,-1.3) rectangle (1,1.3);

    \node at (0,.7) (top) [rectangle,draw] {$  \JW{2n-2} $};
    \node at (0,-.7) (bottom) [ellipse,draw]{$ \hspace{.15 cm} S \hspace{.15cm} $};
    \filldraw (bottom.180) circle (.5mm);    

    \draw (top.-60)--(bottom.60);
    \draw (top.-145)--(bottom.145);
    \node at (-.1,0) {...};
    \draw (top.-35) arc (-180:0:.2cm) -- ++(90:1cm);
    \draw (bottom.35) arc (180:0:.2cm) -- ++(-90:1cm);

    \draw (top.35) -- ++(90:1cm);
    \draw (top.60) -- ++(90:1cm);
    \node at (-.1,1.2) {...};
    \draw (top.145) -- ++(90:1cm);

    \draw (bottom.-35) -- ++(-90:1cm);
    \draw (bottom.-60) -- ++(-90:1cm);
    \node at (-.1,-1.2) {...};
    \draw (bottom.-145) -- ++(-90:1cm);
\end{tikzpicture}
+
\begin{tikzpicture}[baseline]
    \clip (-.8,-1.3) rectangle (1,1.3);

    \node at (0,.7) (top) [ellipse,draw] {$ \hspace{.15 cm} S \hspace{.15cm} $};
    \node at (0,-.7) (bottom) [rectangle,draw] {$  \JW{2n-2} $};
    \filldraw (top.180) circle (.5mm);

    \draw (top.-60)--(bottom.60);
    \draw (top.-145)--(bottom.145);
    \node at (-.1,0) {...};
    \draw (top.-35) arc (-180:0:.2cm) -- ++(90:1cm);
    \draw (bottom.35) arc (180:0:.2cm) -- ++(-90:1cm);

    \draw (top.35) -- ++(90:1cm);
    \draw (top.60) -- ++(90:1cm);
    \node at (-.1,1.2) {...};
    \draw (top.145) -- ++(90:1cm);

    \draw (bottom.-35) -- ++(-90:1cm);
    \draw (bottom.-60) -- ++(-90:1cm);
    \node at (-.1,-1.2) {...};
    \draw (bottom.-145) -- ++(-90:1cm);
\end{tikzpicture}
+\frac{\qi{2n-1}}{\qi{2n}}
\begin{tikzpicture}[baseline]
    \clip (-1,-1.3) rectangle (1,1.3);
    \node at (0,0) (box) [rectangle,draw] {$\hspace{.1cm} \JW{2n-1} \hspace{.1cm} $};

    \draw (box.-50)-- ++(-90:1cm);
    \draw (box.-90) -- ++(-90:1cm);
    \draw (box.-150)--++(-90:1cm);
    \node at (-.25,-.6) {...};
    \draw (box.-30) -- ++(-90:1cm);

    \draw (box.50)-- ++(90:1cm);
    \draw (box.90) -- ++(90:1cm);
    \draw (box.150)--++(90:1cm);
    \node at (-.25,.6) {...};
    \draw (box.30) -- ++(90:1cm);
\end{tikzpicture}
\right)

%% file: diagrams/tikz/PEnPcancel.tex
&=\frac{ \qi{2}}{4} \cdot
\left(
\frac{2}{\qi{2}}
\begin{tikzpicture}[baseline]
    \clip (-1,-1.3) rectangle (1,1.3);
    \node at (0,0) (box) [rectangle,draw] {$\JW{2n-2}$};

    \draw (box.-145)--++(-90:1cm);
    \node at (-.1,-.6) {...};
    \draw (box.-35) -- ++(-90:1cm);
    \draw (box.-60) -- ++(-90:1cm);

    \draw (box.35) -- ++(90:1cm);
    \draw (box.60) -- ++(90:1cm);
    \node at (-.1,.6) {...};
    \draw (box.145) -- ++(90:1cm);

    \draw (.9,1.5)--(.9,-1.5);
\end{tikzpicture}
+
\begin{tikzpicture}[baseline]
    \clip (-.8,-1.3) rectangle (1,1.3);

    \node at (0,.7) (top) [rectangle,draw] {$  \JW{2n-2} $};
    \node at (0,-.7) (bottom) [ellipse,draw]{$ \hspace{.15 cm} S \hspace{.15cm} $};
    \filldraw (bottom.180) circle (.5mm);    

    \draw (top.-60)--(bottom.60);
    \draw (top.-145)--(bottom.145);
    \node at (-.1,0) {...};
    \draw (top.-35) arc (-180:0:.2cm) -- ++(90:1cm);
    \draw (bottom.35) arc (180:0:.2cm) -- ++(-90:1cm);

    \draw (top.35) -- ++(90:1cm);
    \draw (top.60) -- ++(90:1cm);
    \node at (-.1,1.2) {...};
    \draw (top.145) -- ++(90:1cm);

    \draw (bottom.-35) -- ++(-90:1cm);
    \draw (bottom.-60) -- ++(-90:1cm);
    \node at (-.1,-1.2) {...};
    \draw (bottom.-145) -- ++(-90:1cm);
\end{tikzpicture}
+
\begin{tikzpicture}[baseline]
    \clip (-.8,-1.3) rectangle (1,1.3);

    \node at (0,.7) (top) [ellipse,draw] {$ \hspace{.15 cm} S \hspace{.15cm} $};
    \node at (0,-.7) (bottom) [rectangle,draw] {$  \JW{2n-2} $};
    \filldraw (top.180) circle (.5mm);

    \draw (top.-60)--(bottom.60);
    \draw (top.-145)--(bottom.145);
    \node at (-.1,0) {...};
    \draw (top.-35) arc (-180:0:.2cm) -- ++(90:1cm);
    \draw (bottom.35) arc (180:0:.2cm) -- ++(-90:1cm);

    \draw (top.35) -- ++(90:1cm);
    \draw (top.60) -- ++(90:1cm);
    \node at (-.1,1.2) {...};
    \draw (top.145) -- ++(90:1cm);

    \draw (bottom.-35) -- ++(-90:1cm);
    \draw (bottom.-60) -- ++(-90:1cm);
    \node at (-.1,-1.2) {...};
    \draw (bottom.-145) -- ++(-90:1cm);
\end{tikzpicture}
\right)

%% file: diagrams/tikz/JWEnS.tex
\begin{tikzpicture}[baseline]
    \clip (-.8,-1.3) rectangle (1,1.3);

    \node at (0,.7) (top) [rectangle,draw] {$  \JW{2n-2} $};
    \node at (0,-.7) (bottom) [ellipse,draw]{$ \hspace{.15 cm} S \hspace{.15cm} $};
    \filldraw (bottom.180) circle (.5mm);    

    \draw (top.-60)--(bottom.60);
    \draw (top.-145)--(bottom.145);
    \node at (-.1,0) {...};
    \draw (top.-35) arc (-180:0:.2cm) -- ++(90:1cm);
    \draw (bottom.35) arc (180:0:.2cm) -- ++(-90:1cm);

    \draw (top.35) -- ++(90:1cm);
    \draw (top.60) -- ++(90:1cm);
    \node at (-.1,1.2) {...};
    \draw (top.145) -- ++(90:1cm);

    \draw (bottom.-35) -- ++(-90:1cm);
    \draw (bottom.-60) -- ++(-90:1cm);
    \node at (-.1,-1.2) {...};
    \draw (bottom.-145) -- ++(-90:1cm);
\end{tikzpicture}
+
\begin{tikzpicture}[baseline]
    \clip (-.8,-1.3) rectangle (1,1.3);

    \node at (0,.7) (top) [ellipse,draw] {$ \hspace{.15 cm} S \hspace{.15cm} $};
    \node at (0,-.7) (bottom) [rectangle,draw] {$  \JW{2n-2} $};
    \filldraw (top.180) circle (.5mm);

    \draw (top.-60)--(bottom.60);
    \draw (top.-145)--(bottom.145);
    \node at (-.1,0) {...};
    \draw (top.-35) arc (-180:0:.2cm) -- ++(90:1cm);
    \draw (bottom.35) arc (180:0:.2cm) -- ++(-90:1cm);

    \draw (top.35) -- ++(90:1cm);
    \draw (top.60) -- ++(90:1cm);
    \node at (-.1,1.2) {...};
    \draw (top.145) -- ++(90:1cm);

    \draw (bottom.-35) -- ++(-90:1cm);
    \draw (bottom.-60) -- ++(-90:1cm);
    \node at (-.1,-1.2) {...};
    \draw (bottom.-145) -- ++(-90:1cm);
\end{tikzpicture}
=
\frac{2}{\qi{2}}
\begin{tikzpicture}[baseline]
    \clip (-.7,-1.3) rectangle (1,1.3);
    \node at (0,0) (box) [ellipse,draw] {$ \hspace{.15 cm} S \hspace{.15cm} $};
    \filldraw (box.180) circle (.5mm);

    \draw (box.-145)--++(-90:1cm);
    \node at (box.-90) [below] {...};
    \draw (box.-35) -- ++(-90:1cm);

    \draw (box.35) -- ++(90:1cm);
    \node at (box.90) [above] {...};
    \draw (box.145) -- ++(90:1cm);

    \draw (.9,1.5)--(.9,-1.5);
\end{tikzpicture}

%% file: diagrams/tikz/JWEnStopcap.tex
\begin{tikzpicture}[baseline]
    \clip (-.8,-1.3) rectangle (1,1.3);

    \node at (0,.7) (top) [rectangle,draw] {$  \JW{2n-2} $};
    \node at (0,-.7) (bottom) [ellipse,draw]{$ \hspace{.15 cm} S \hspace{.15cm} $};
    \filldraw (bottom.180) circle (.5mm);    

    \draw (top.-60)--(bottom.60);
    \draw (top.-145)--(bottom.145);
    \node at (-.1,0) {...};
    \draw (top.-35) arc (-180:0:.2cm) -- ++(90:.7cm) arc (0:180:.2cm) --(top.35);
    \draw (bottom.35) arc (180:0:.2cm) -- ++(-90:1cm);

    \draw (top.60) -- ++(90:1cm);
    \node at (-.1,1.2) {...};
    \draw (top.145) -- ++(90:1cm);

    \draw (bottom.-35) -- ++(-90:1cm);
    \draw (bottom.-60) -- ++(-90:1cm);
    \node at (-.1,-1.2) {...};
    \draw (bottom.-145) -- ++(-90:1cm);
\end{tikzpicture}
+
\begin{tikzpicture}[baseline]
    \clip (-.8,-1.3) rectangle (1,1.3);

    \node at (0,.7) (top) [ellipse,draw] {$ \hspace{.15 cm} S \hspace{.15cm} $};
    \node at (0,-.7) (bottom) [rectangle,draw] {$  \JW{2n-2} $};
    \filldraw (top.180) circle (.5mm);

    \draw (top.-60)--(bottom.60);
    \draw (top.-145)--(bottom.145);
    \node at (-.1,0) {...};
    \draw (top.-35) arc (-180:0:.2cm) -- ++(90:.6cm) arc (0:180:.2cm) --(top.35);
    \draw (bottom.35) arc (180:0:.2cm) -- ++(-90:1cm);

    \draw (top.60) -- ++(90:1cm);
    \node at (-.1,1.2) {...};
    \draw (top.145) -- ++(90:1cm);

    \draw (bottom.-35) -- ++(-90:1cm);
    \draw (bottom.-60) -- ++(-90:1cm);
    \node at (-.1,-1.2) {...};
    \draw (bottom.-145) -- ++(-90:1cm);
\end{tikzpicture}
& = \frac{2}{\qi{2}}
\begin{tikzpicture}[baseline]
    \clip (-.7,-1.3) rectangle (1,1.3);
    \node at (0,0) (box) [ellipse,draw] {$ \hspace{.15 cm} S \hspace{.15cm} $};
    \filldraw (box.180) circle (.5mm);

    \draw (box.-145)--++(-90:1cm);
    \node at (-.1,-.6) {...};
    \draw (box.-35) -- ++(-90:1cm);
    \draw (box.-60) -- ++(-90:1cm);

    \draw (box.35) arc (180:0:.2cm) -- ++(-90:2cm);
    \draw (box.60) -- ++(90:1cm);
    \node at (-.1,.6) {...};
    \draw (box.145) -- ++(90:1cm);
\end{tikzpicture}

%% file: diagrams/tikz/JWEnSbottomcap.tex
\begin{tikzpicture}[baseline]
    \clip (-.8,-1.3) rectangle (1,1.3);

    \node at (0,.7) (top) [rectangle,draw] {$  \JW{2n-2} $};
    \node at (0,-.7) (bottom) [ellipse,draw]{$ \hspace{.15 cm} S \hspace{.15cm} $};
    \filldraw (bottom.180) circle (.5mm);    

    \draw (top.-60)--(bottom.60);
    \draw (top.-145)--(bottom.145);
    \node at (-.1,0) {...};
    \draw (top.-35) arc (-180:0:.2cm) -- ++(90:1cm);
    \draw (bottom.-35) arc (-180:0:.2cm) -- ++(90:.6cm) arc (0:180:.2cm) --(bottom.35);

    \draw (top.60) -- ++(90:1cm);
    \node at (-.1,1.2) {...};
    \draw (top.145) -- ++(90:1cm);

    \draw (top.35) -- ++(90:1cm);
    \draw (bottom.-60) -- ++(-90:1cm);
    \node at (-.1,-1.2) {...};
    \draw (bottom.-145) -- ++(-90:1cm);
\end{tikzpicture}
+
\begin{tikzpicture}[baseline]
    \clip (-.8,-1.3) rectangle (1,1.3);

    \node at (0,.7) (top) [ellipse,draw] {$ \hspace{.15 cm} S \hspace{.15cm} $};
    \node at (0,-.7) (bottom) [rectangle,draw] {$  \JW{2n-2} $};
    \filldraw (top.180) circle (.5mm);

    \draw (top.-60)--(bottom.60);
    \draw (top.-145)--(bottom.145);
    \node at (-.1,0) {...};
    \draw (top.-35) arc (-180:0:.2cm) -- ++(90:1cm);
    \draw (bottom.35) arc (180:0:.2cm) -- ++(-90:.7cm) arc (0:-180:.2cm) --(bottom.-35);

    \draw (top.60) -- ++(90:1cm);
    \node at (-.1,1.2) {...};
    \draw (top.145) -- ++(90:1cm);

    \draw (top.35) -- ++(90:1cm);
    \draw (bottom.-60) -- ++(-90:1cm);
    \node at (-.1,-1.2) {...};
    \draw (bottom.-145) -- ++(-90:1cm);
\end{tikzpicture}
& = \frac{2}{\qi{2}}
\begin{tikzpicture}[baseline]
    \clip (-.7,-1.3) rectangle (1,1.3);
    \node at (0,0) (box) [ellipse,draw] {$ \hspace{.15 cm} S \hspace{.15cm} $};
    \filldraw (box.180) circle (.5mm);

    \draw (box.-145)--++(-90:1cm);
    \node at (-.1,-.6) {...};
    \draw (box.35) -- ++(90:1cm);
    \draw (box.-60) -- ++(-90:1cm);

    \draw (box.-35) arc (-180:0:.2cm) -- ++(90:2cm);
    \draw (box.60) -- ++(90:1cm);
    \node at (-.1,.6) {...};
    \draw (box.145) -- ++(90:1cm);
\end{tikzpicture}

%% file: diagrams/tikz/JWS.tex
\begin{tikzpicture}[baseline]
    \clip (-.8,-1.3) rectangle (1,1.3);

    \node at (0,.7) (top) [rectangle,draw] {$  \JW{2n-2} $};
    \node at (0,-.7) (bottom) [ellipse,draw]{$ \hspace{.15 cm} S \hspace{.15cm} $};
    \filldraw (bottom.180) circle (.5mm);    

    \draw (top.-60)--(bottom.60);
    \draw (top.-145)--(bottom.145);
    \node at (-.1,0) {...};
    \draw (top.-35) -- (bottom.35);

    \draw (top.35) -- ++(90:1cm);
    \draw (top.60) -- ++(90:1cm);
    \node at (-.1,1.2) {...};
    \draw (top.145) -- ++(90:1cm);

    \draw (bottom.-35) -- ++(-90:1cm);
    \draw (bottom.-60) -- ++(-90:1cm);
    \node at (-.1,-1.2) {...};
    \draw (bottom.-145) -- ++(-90:1cm);
\end{tikzpicture}
+
\begin{tikzpicture}[baseline]
    \clip (-.8,-1.3) rectangle (1,1.3);

    \node at (0,.7) (top) [ellipse,draw] {$ \hspace{.15 cm} S \hspace{.15cm} $};
    \node at (0,-.7) (bottom) [rectangle,draw] {$  \JW{2n-2} $};
    \filldraw (top.180) circle (.5mm);

    \draw (top.-60)--(bottom.60);
    \draw (top.-145)--(bottom.145);
    \node at (-.1,0) {...};
    \draw (top.-35) -- (bottom.35);

    \draw (top.35) -- ++(90:1cm);
    \draw (top.60) -- ++(90:1cm);
    \node at (-.1,1.2) {...};
    \draw (top.145) -- ++(90:1cm);

    \draw (bottom.-35) -- ++(-90:1cm);
    \draw (bottom.-60) -- ++(-90:1cm);
    \node at (-.1,-1.2) {...};
    \draw (bottom.-145) -- ++(-90:1cm);
\end{tikzpicture}
& = \frac{2}{\qi{2}}
\begin{tikzpicture}[baseline]
    \clip (-.7,-1.3) rectangle (1.5,1.3);
    \node at (0,0) (box) [ellipse,draw] {$ \hspace{.15 cm} S \hspace{.15cm} $};
    \filldraw (box.180) circle (.5mm);

    \draw (box.-145)--++(-90:1cm);
    \node at (box.-90) [below] {...};
    \draw (box.-35) -- ++(-90:1cm);

    \draw (box.35) -- ++(90:1cm);
    \node at (box.90) [above] {...};
    \draw (box.145) -- ++(90:1cm);

    \draw (.9,.5)--(.9,-.5) arc (-180:0:.2cm) -- ++(90:1cm) arc (0:180:.2cm);
\end{tikzpicture}
=
2
\begin{tikzpicture}[baseline]
    \clip (-.7,-1.3) rectangle (.7,1.3);
    \node at (0,0) (box) [ellipse,draw] {$ \hspace{.15 cm} S \hspace{.15cm} $};
    \filldraw (box.180) circle (.5mm);

    \draw (box.-145)--++(-90:1cm);
    \node at (box.-90) [below] {...};
    \draw (box.-35) -- ++(-90:1cm);

    \draw (box.35) -- ++(90:1cm);
    \node at (box.90) [above] {...};
    \draw (box.145) -- ++(90:1cm);
\end{tikzpicture}

%% file: text/basis.tex
\subsection{A basis for the planar algebra}%
\label{sec:basis}%
In this section we present an explicit basis for the
planar algebra, and use this to show that the generators and relations
presentation from Definition \ref{def:pa} really does result in a
positive definite planar algebra.

Each vector space $\pa_m$ of the planar algebra also appears as a
$\operatorname{Hom}$ space of the corresponding tensor category of
projections, specifically as $\Hom{}{\id}{X^{\tensor m}}$. We'll use a standard approach for describing bases for
semisimple tensor categories, based on tree-diagrams.

For each triple of self-adjoint minimal projections $p,q,r$, we need to fix an
orthogonal basis for $\Hom{}{\id}{p \tensor q \tensor r}$. Call these bases $\{v_\lambda\}_{\lambda \in \mathcal{B}(p,q,r)}$.
If we take the adjoint of $v_\lambda \in \Hom{}{\id}{p \tensor q \tensor r}$, we get
$v_\lambda^* \in \Hom{}{p \tensor q \tensor r}{\id}$.

In fact, we'll only need to do this when one of the the three projections $p, q$ and $r$ is just $X$. In these cases, we've already implicitly described
the $\operatorname{Hom}$ spaces in Lemmas \ref{LemmasB1}, \ref{LemmasB2} and \ref{LemmasB3}.

We can now interpret certain planar trivalent graphs as notations for elements
of the planar algebra. The graphs have oriented edges labelled by projections, but where we allow
reversing the orientation and replacing the projection with its dual.\footnote{We'll leave off orientations on edges labelled by $X$, since it's self dual. In fact, all the minimal projections are self dual, except for $P$ and $Q$ when $n$ is even, in which case $\overline{P} = Q$ and $\overline{Q}=P$.} The graphs have vertices labelled by elements of the sets $\mathcal{B}(p,q,r)$
described above (where $p, q$ and $r$ are the projections on the edges
\emph{leaving} the vertex).  If $\card{\mathcal{B}(p,q,r)}=1$ we may leave off the label at that vertex.

To produce an element of the planar algebra from such a graph, we simply
replace each edge labelled by a projection $p$ in $\pa_{2m}$ with $m$
parallel strands, with the projection $p$ drawn across them, and each
trivalent vertex labelled by $\lambda$ with the element $v_\lambda \in
\Hom{}{\id}{p \tensor q \tensor r}$.

As a first example,
\begin{defn}
We call the norm of the element $v_\lambda$, with $\lambda \in
\mathcal{B}(p,q,r)$, the theta-symbol:

\begin{equation*}
\theta(p,q,r;\lambda) := \mathfig{0.2}{basis/theta}
\end{equation*}
\end{defn}

\begin{defn}
Fix a list of minimal projections $(p_i)_{0\leq i\leq k+1}$, called the
boundary. A tree diagram for this boundary is a trivalent graph of the form:
\begin{equation*}
\mathfig{0.75}{basis/tree}
\end{equation*}

It is labelled by
\begin{itemize}
\item another list of minimal projections $(q_i)_{1\leq i \leq k-1}$
such that $q_i$ is a summand of $q_{i-1} \tensor p_i$ for each $1\leq i
\leq k$, or equivalently that $\mathcal{B}(\overline{q_i}, q_{i-1}, p_i) \neq
\eset$ (here we make the identifications $q_0 = p_0$ and $q_k =
\overline{p_{k+1}}$),
\item and for each $1\leq i\leq k$, a choice of orthogonal basis vector
$v_{\lambda_i}$, with $\lambda_i \in \mathcal{B}(\overline{q_i}, q_{i-1}, p_i)$.
\end{itemize}
\end{defn}

\begin{thm}
\label{thm:identity-tree}
The $k$-strand identity can be written as a sum of tree diagrams.
(We'll assume there are no multiple edges in the principal graph for the exposition here; otherwise, we need to remember labels at vertices.)
Let $\Gamma_{k-1}$ be the set of length $k-1$ paths on the principal graph starting at $X$. (Thus if $\gamma \in \Gamma_{k-1}$, $\gamma_0 = X$ and the endpoint of the path is $\gamma_{k-1}$.)
Then
\begin{equation*}
\mathfig{0.18}{basis/identity} = \sum_{\gamma \in \Gamma_{k-1}} \prod_{i=0}^{k-2} \frac{\tr{\gamma_{i+1}}}{\theta(\overline{\gamma_{i}},\gamma_{i+1},X))} \mathfig{0.18}{basis/identity-term}
\end{equation*}
\end{thm}
\begin{proof}
We induct on $k$.

When $k=1$, the result is trivially true; the only path in $\Gamma_0$ is the constant path, with $\gamma(0)=X$, and there's no coefficient.

To prove the result for $k+1$, we replace the first $k$ strands on the left, obtaining
\begin{equation*}
\mathfig{0.18}{basis/identity_1} = \sum_{\gamma \in P_k} \prod_{i=0}^{k-2} \frac{\tr{\gamma_{i+1}}}{\theta(\overline{\gamma_{i}},\gamma_{i+1},X)} \mathfig{0.18}{basis/identity-term_1}
\end{equation*}
and then use the identity
\begin{equation*}
\mathfig{0.075}{basis/gamma_1} = \sum_{\substack{\text{$\gamma_k$ adjacent}\\\text{to $\gamma_{k-1}$}}} \frac{\tr{\gamma_k}}{\theta(\overline{\gamma_{k-1}},\gamma_k,X)} \mathfig{0.075}{basis/gamma-term}
\end{equation*}
(which certainly holds with some coefficients, by the definition of the principal graph, and with these particular coefficients by multiplying in turn both sides by each of the terms on the right)
to obtain the desired result.
\end{proof}

\begin{thm}
The tree diagrams with boundary labelled entirely by $X$ give a positive orthogonal basis for the
invariant space $\Hom{}{\id}{X^{\tensor n}}$.
\end{thm}
\begin{rem}
Actually, the tree diagrams with boundary $(p_i)$ give an orthogonal basis for the invariant space
$\Hom{}{\id}{\Tensor_i p_i}$, but we won't prove that here. We'd need to exhibit explicit bases for all the triple invariant spaces in order
to check positivity, and a slightly stronger version of Theorem \ref{thm:identity-tree}.
\end{rem}
\begin{proof}
To see that the tree diagrams are all orthogonal is just part of the
standard machinery of semisimple tensor categories --- make repeated use
of the formulas
\begin{align*}
\mathfig{0.08}{basis/tr-p} & = \tr{p} \\
\mathfig{0.12}{basis/open-theta} & = \delta_{p=q} \delta_{\mu=\lambda^*} \frac{\theta(\overline{p},r,s;\lambda)}{\tr{p}} \mathfig{0.032}{basis/strand}
\end{align*}
where $\lambda \in \mathcal{B}(\overline{p},r,s)$ and $\mu \in \mathcal{B}(\overline{s},\overline{r},q)$.
(In fact, this proves that the tree diagrams are orthogonal for arbitrary boundaries).

The norm of a tree diagram is a ratio of theta symbols and traces of
projections.
The trace of $f^{(k)}$ is $\qi{k+1}$, and $\tr{P} = \tr{Q} = \frac{\qi{2n-2}}{2}$, and these quantities are all positive at our value of $q$. Further, the theta symbols with one edge labelled by $X$ are all easy to calculate (recall the relevant one-dimensional bases for the $\operatorname{Hom}$ spaces were described in Lemmas \ref{LemmasB1}, \ref{LemmasB2} and \ref{LemmasB3}), and
in fact are just equal to traces of these same projections:
\begin{align*}
\theta(f^{(k-1)},f^{(k)},X) & = \tr{f^{(k)}} \\
\theta(f^{(2n-3)},P,X) & = \tr{P} \\
\theta(f^{(2n-3)},Q,X) & = \tr{Q} \\
\end{align*}
Since these are all positive, the norms of tree diagrams are positive.

To see that the tree diagrams span, we make use of Theorem \ref{thm:identity-tree}, and Lemma \ref{lem:no-homs}.
Take an arbitrary open diagram $D$ with $k$ boundary points, and write it as $D \cdot \id_k$. Apply Theorem \ref{thm:identity-tree} to $\id_k$, and
observe that all terms indexed by paths not ending at $f^{(0)}$ are zero, by Lemma \ref{lem:no-homs} (here we think of $D$ as having an extra boundary point labeled by $f^{(0)}$, so we get a map from $f^{(0)}$ to the endpoint of the path).
In the remaining terms, we have the disjoint union (after erasing the innermost edge labeled by $f^{(0)}$) of a closed diagram
and a tree diagram. Since all closed diagrams can be evaluated, by Corollary \ref{cor:evaluation}, we see we have rewritten an arbitrary diagram
as a linear combination of tree diagrams.
\end{proof}

\begin{cor}
The planar algebra given by generators and relations in Definition
\ref{def:pa} is positive definite.
\end{cor}

Therefore, $\pa$ is indeed the subfactor planar algebra with principal graph $D_{2n}$.

%% file: text/appendix.tex
\appendix

\section{A brief note on $T_n$, a related planar algebra}
\label{appendix}
In this section we briefly describe modifications of the skein relations
for $D_{2n}$ which give rise to the planar algebras $T_n$. The planar algebras $T_n$ have appeared previously in \cite{MR1936496, MR2046203, MR1333750}.
They are \emph{unshaded} subfactor planar algebras in the sense
we've described in \ref{subsec:pa}, but they are not \emph{shaded} subfactor planar algebras (the more usual sense).

The most direct construction of the $T_n$ planar algebra is to interpret the single strand as $f^{(2n-2)}$ in the Temperley-Lieb planar algebra $A_{2n}$,
allowing arbitrary Temperley-Lieb diagrams with $(2n-2) m$ boundary points in the $m$-boxes. (Another way to say this, in the langauge of tensor categories with a distinguished tensor generator,
is to take the even subcategory of $A_{2n}$, thought of as generated by $f^{(2n-2)}$.) This certainly ensures that $T_n$ exists; below
we give a presentation by generators and relations.

We consider a skein theory with a $(k=2n+1)$ strand generator (allowing
in this appendix boxes with odd numbers of boundary points), at the
special value $q=e^{\frac{i\pi}{k+2}}$, and relations analogous to those of Definition \ref{def:pa}:
\begin{enumerate}

\item\label{delta-T} a closed loop is equal to $2 \cos(\frac{\pi}{k+2})$,

\item\label{rotateS-T}
\input{diagrams/tikz/rotateS.tex}

\item\label{capS-T}
\input{diagrams/tikz/capS-T.tex}

\item\label{twoS-T}
\input{diagrams/tikz/twoS-T.tex}
\end{enumerate}

A calculation analogous to that of Theorem \ref{thm:passacrossS} shows that we have the relations
\begin{equation*}
\input{diagrams/tikz/pullstringoverS-T.tex}
\qquad \text{and} \qquad \input{diagrams/tikz/pullstringunderS-T.tex}.
\end{equation*}
where $Z^- = Z^+ = (-1)^{\frac{k+1}{2}}$. (Recall that in the $D_{2n}$ case discussed in the body of the paper we had $Z^\pm = \pm 1$.)
These relations allow us to repeat the arguments showing that closed diagrams can be evaluated, and that the planar algebra is spherical.
When $k \equiv 3 \pmod{4}$ and $Z^\pm = +1$, the planar algebra $T_n$ is braided.

When $k \equiv 1 \pmod{4}$ and $Z^\pm = -1$, one can replace the usual crossing in Temperley-Lieb with minus itself; this is still a braiding on
Temperley-Lieb. One then has instead $Z^\pm = +1$, and so the entire planar algebra is then honestly braided.
Notice that $T_n$ is related to $A_{2n}$ in two \emph{different} ways. First, $T_n$ contains $A_{2n}$ as a subplanar algebra (simply because any planar algebra at a special value
of $\qi{2}$ contains the corresponding Temperley-Lieb planar algebra). Second, $T_n$ is actually the even part of $A_{2n}$, with an
unusual choice of generator (see above). The first gives a candidate braiding -- as we've seen it's only an `almost braiding' when $k \equiv 1 \pmod{4}$.
The second automatically gives an honest braiding, and in the $k \equiv 1 \pmod{4}$ case it's the negative of the first one.

Following through the consistency argument of \S \ref{sec:consistency}, mutatis mutandi,
we see that these relations do not collapse the planar algebra to zero. Further, along the lines of \S \ref{sec:category}, we can show that the
tensor category of projections is semisimple, with $\{f^{(0)}, f^{(1)}, \ldots, f^{(\frac{k-1}{2})}\}$ forming a complete orthogonal set of minimal projections.
The element $S$ in the planar algebra gives rise to isomorphisms $f^{(i)} \iso f^{(k-i)}$ for $i= 0, \ldots, \frac{k-1}{2}$. Further, the principal graph
is $T_{\frac{k-1}{2}}$, the tadpole graph:
\begin{equation*}
\mathfig{0.5}{graphs/Tk}.
\end{equation*}

%% file: diagrams/tikz/capS-T.tex
$
\begin{tikzpicture}[baseline]
 \node (S) [circle,draw] {$S$};
  \filldraw (S.157) circle (.5mm);
 \draw (S.135) .. controls +(135:3mm) and +(90:3mm) .. (S.90);
 \draw (S.45) -- +(45:3mm);
 \draw (S.0) -- +(0:3mm);
 \draw (S.-45) node [below] {$\cdot$};
  \draw (S.-90) node [below] {$\cdot$};
 \draw (S.-135) node [below] {$\cdot$};
 \draw (S.180) -- +(180:3mm);
\end{tikzpicture}
= 0
$

%% file: diagrams/tikz/twoS-T.tex
$
\begin{tikzpicture}[baseline]
 \node at (0,.7) (S) [circle,draw] {$S$};
     \filldraw (S.180) circle (.5mm);

 \node at (0,-.7) (S') [circle,draw] {$S$};
    \filldraw (S'.180) circle (.5mm);

 \draw (S.45) -- +(45:4mm);
 \draw (S.90) node [above] {$\ldots$};
 \draw (S.135) -- +(135:4mm);
 \draw (S'.225) -- +(-135:4mm);
 \draw (S'.270) node [below] {$\ldots$};
 \draw (S'.-45) -- +(-45:4mm);
\end{tikzpicture}
=\qi{n+1} \cdot
\begin{tikzpicture}[baseline]
\node (0,0) (JW) [rectangle,draw] {$\quad \JW{k} \quad§·$};
\draw (JW.180);

\draw (JW.35) -- +(90:11mm);
\draw (JW.90) node [above] {$\cdots$};
\draw (JW.145) -- +(90:11mm);

\draw (JW.-35) -- +(-90:11mm);
\draw (JW.-90) node [below] {$\cdots$};
\draw (JW.-145) -- +(-90:11mm);
\end{tikzpicture}
$

%% file: diagrams/tikz/pullstringoverS-T.tex
\begin{tikzpicture}[baseline]
    \clip (-1.1,-1.4) rectangle (1.1,1.4);

    \node (S) at (0,0) [circle, draw] {$S$};
     \filldraw (S.180) circle (.5mm);

    \draw[rounded corners=2mm] (-1,2) -- (-1,-1) -- (1,-1) -- (1,2);

    \draw (S.135) -- ++(90:15mm);
    \draw (S.90) node [above] {...};
    \draw (S.45) -- ++(90:15mm);
\end{tikzpicture}
=Z_k^+
\begin{tikzpicture}[baseline]
    \clip (-1.1,-1.4) rectangle (1.1,1.4);

    \node (S) at (0,0) [circle, draw] {$S$};
     \filldraw (S.180) circle (.5mm);
    \node (x1) at (S.135 |- -1,1){};
    \node (x2) at  (S.45 |- -1,1){} ;

    \draw[rounded corners=2mm] (-1,1.5) --  (-1,1) -- (1,1) -- (1,1.5);

    \draw (S.135) -- (x1.-90);
    \draw (x1.90) -- ++(90:5mm);
    \draw (S.90) node [above] {...};
    \draw (x2.90) -- ++(90:5mm);
    \draw (S.45) -- (x2.-90);
\end{tikzpicture}

%% file: diagrams/tikz/pullstringunderS-T.tex
\begin{tikzpicture}[baseline]
    \clip (-1.1,-1.4) rectangle (1.1,1.4);

    \node (S) at (0,0) [circle, draw] {$ S$};
     \filldraw (S.180) circle (.5mm);

    \draw[rounded corners=2mm] (-1,2) -- (-1,-1) -- (1,-1) -- (1,2);

    \draw (S.135) -- ++(90:15mm);
    \draw (S.90) node [above] {...};
    \draw (S.45) -- ++(90:15mm);
\end{tikzpicture}
=Z_k^-
\begin{tikzpicture}[baseline]
    \clip (-1.1,-1.4) rectangle (1.1,1.4);

    \node (S) at (0,0) [circle, draw] {$S$};
     \filldraw (S.180) circle (.5mm);

    \node (x1) at (S.135 |- -1,1){};
    \node (x2) at  (S.45 |- -1,1){} ;

    \draw[rounded corners=2mm](-1,1.5) -- (-1,1) -- (x1.180);
    \draw (x1.0) -- (x2.180);
    \draw[rounded corners=2mm] (x2.0) -- (1,1)--(1,1.5);

    \draw (S.135) -- (x1.90);
    \draw (x1.90) -- ++(90:5mm);
    \draw (S.90) node [above] {...};
    \draw (x2.90) -- ++(90:5mm);
    \draw (S.45) -- (x2.90);
\end{tikzpicture}